\documentclass[notitlepage,leqno,10pt]{article}

\newenvironment{proof}{\noindent{\sc Proof}.\enspace}{\rule{2mm}{2mm}\medskip}

\usepackage{amsmath}
\usepackage{amssymb}

\newtheorem{theorem}{Theorem}[section]
\newtheorem{lemma}[theorem]{Lemma}

\newtheorem{remark}[theorem]{Remark}
\newtheorem{example}[theorem]{Example}

\usepackage{amsmath}
\usepackage{amssymb}

\newcommand{\loc}{{loc}}
\newcommand{\capacity}{{\mathcal C}}
\newcommand{\capinRn}{{\mathbf{C}}}

\newcommand{\diam}{{ \rm diam}}
\newcommand{\supp}{{ \rm supp}}
\newcommand{\dist}{{ \rm dist}}
\newcommand{\kruzhok}{ \ensuremath{\mbox{\,}_\circ}}
\newcommand{\gstand}{ \ensuremath{\overset{\kruzhok}{g}}}
\newcommand{\geucl}{ \ensuremath{g_E}}
\newcommand{\Hmeas}{{ \mathcal{H}}}
\newcommand{\metricexp}{\ensuremath{{2/(n-2)}}}

\title{Thinness for Scalar-Negative Singular Yamabe Metrics}

\author{Denis A. Labutin}

\date{}

\begin{document}
\sloppy

\maketitle

\begin{center}

{\small
Department of Mathematics, University of California,
Santa Barbara, CA 93106, USA}

\end{center}

\footnotetext[1]{E-mail address: 
labutin@math.ucsb.edu}

\

\

\noindent
{\sc abstract}.
This paper deals with the conformal 
deformation of the standard metric 
in a domain on the  sphere
to a complete metric with the 
constant  scalar curvature.
The problem of description of  domains 
allowing such deformation
originates in the works of  
Loewner and Nirenberg, 
and Schoen and Yau
concerned with the locally conformally flat manifolds.
The goal of this work is to
apply  ideas  from the nonlinear 
potential theory to the problem.
They allow,
in particular,
to solve the problem
in the case of the 
constant negative scalar curvature.

\begin{center}

\bigskip\bigskip

\centerline{\bf AMS subject classification: 
35J60,  53A30, 58J05
}

\end{center}

%
%
%
%
%
%
%
%
%
%
%
%
%
%
%
%

\section{Introduction}
\label{sec1}

\setcounter{equation}{0}

\subsection{Singular Yamabe problem}
\label{sec1.1}

Yamabe problem 
\cite{Yam}
was to prove that for any  compact Riemannian
manifold
$(M, g)$
of dimension
$n\geq 3$
we can find a 
metric conformal to 
$g$
with a constant scalar curvature
by solving a certain variational problem.
This  was proved   in three subsequent 
contributions by
Trudinger 
\cite{Tru},
Aubin 
\cite{Au1},
and Schoen
\cite{Schoen_Yamabe},
see also 
\cite{Lee_Parker},
\cite{Schoen_Yau_book}
for an exposition.
Later other proofs were given,
see
\cite{Aubin_book}
for a survey.
To a large extent,
it is the Yamabe problem that stimulated
the development of the modern geometric analysis.
An intensive work was done on Yamabe problem for manifolds with  boundary.
In this case seeks a conformal metric
with constant curvatures in the interrior and on the boundary, cf. e.g.
\cite{Escobar_92_1},
\cite{Escobar_92_2},
\cite{Escobar_03}.
There were also generalisations to non-Riemannian
settings
\cite{Jerison_Lee},
\cite{Gamara_Yacoub}.

In 1988
Schoen and Yau 
arrived at a different 
(non-variational)
problem
\cite{Schoen_Yau_Invent_88},
\cite{Schoen_ICM}.
Namely, is it possible to characterise domains
on the unit sphere admitting a conformal
deformation of the standard metric to a 
{\it complete}
metric with a constant scalar curvature?
Schoen and Yau were led
to this problem
by their research on geometry and topology
of locally conformally flat manifolds
and were mainly interested in the case of
non-negative curvature.
The case of negative scalar curvature
goes back to an early paper by 
Loewner and Nirenberg
\cite{Loewner_Nirenberg}.
Further motivations for the problem can be 
found  in 
\#36
from  
\cite{Yau_problems},
\cite{Schoen_ICM}
\cite{McOwen_survey}.

The goal of this paper is to introduce methods
of nonlinear potential theory to this problem.
They allow, in particular,
to solve the problem in the negative curvature case.
Our Theorem~\ref{maintheorem}
states that the conformal deformation to a complete
scalar-negative metric is possible in
$\Omega$
if and only if
its complement
is  {\it not thin}.
Thinness, see sec.1.2,  is a basic concept
in  potential theory first introduced by Wiener
in his works on the classical Dirichlet problem.
Developments in nonlinear potential theory
easily allow to relate thinness with geometric properties.
Let us now describe the previous work and our results on the
problem in more details.

Under the conformal
change of metric
${g}=u^{4/(n-2)}\gstand$,
$n\geq 3$,
the scalar curvature changes according  to
the formula
\begin{equation}
\label{scal_curv_change}
R({g})
=
u^{- (n+2)/(n-2)}
\left(
-
\frac{4(n-1)}{n-2}
\Delta u 
+
R(\gstand) u
\right).
\end{equation}
Here 
$
\Delta u= \mathrm{div}\left( \mathrm{grad}\, u\right)
$
is the Laplace-Beltrami  operator
on
$({\bf S}^n, \gstand)$
and
$R(\gstand)=n(n-1)$
is the scalar curvature of 
the standard metric
$\gstand$
induced by the embedding
${\bf S}^n \hookrightarrow {\bf R}^{n+1}$.
Thus
analytically 
for  given
$\Omega\subset{\bf S}^n$
and
$R\in\{-1, 0, 1\}$
one  seeks  
a  smooth  solution
to the following 
problem:
\begin{eqnarray}
\label{singyam}
\frac{4(n-1)}{n-2}
\Delta  u  
-
R(\gstand)
u  
+ 
R
u^{(n+2)/(n-2)}
=0
&
\
{\rm in}
\
&
\Omega
\nonumber
\\
u>0
&
\
{\rm in}
\
&
\Omega
\\
u^{4/(n-2)}\gstand
\
{\rm is \ complete\ metric}
&
\
{\rm in}
\
&
\Omega.
\nonumber
\end{eqnarray}
The case
$R=1$
is regarded as the hardest among the  three.
The present work focuses on the negative curvature case
\begin{equation}
\label{R=-1}
R=-1.
\end{equation}
Survey
\cite{McOwen_survey}
by McOwen
describes the progress on  the problem
and open problems.
Let us explain how our results fit in 
the general picture.
We set
$K={\bf S}^n\setminus\Omega$.

Investigations of the negative curvature case were started
in
1974
by
Loewner and Nirenberg
\cite{Loewner_Nirenberg}.
They proved that 
if the 
problem
(\ref{singyam}),
(\ref{R=-1})
admits a solution then
the complement
of
$\Omega$
must satisfy
\begin{equation*}
\mathcal{H}^{(n-2)/2}(K)=\infty.
\end{equation*}
Here
$\mathcal{H}^\alpha$
denotes the Hausdorff $\alpha$-measure.
Their work together with
Aviles
\cite{Aviles_CPDE}
and Veron
\cite{Veron_JDE}
showed that
if 
$K$
is a smooth submanifold
of
$\mathbf{S}^n$
of the corresponding dimension
$k>(n-2)/2$
then problem
(\ref{singyam}),
(\ref{R=-1})
has a solution.
Mazzeo
\cite{Mazzeo_IUMJ_91}
showed, in particular,   that for such
$K$ 
the solution is unique.
Finn
\cite{Finn_94},
\cite{Finn_99},
\cite{Finn_00}
established the solvability 
under weaker conditions on
$K$.
Namely he required that
it has a structure similar to
(actually, more general than)
Lipschitz submanifold of the corresponding dimension, see also
\cite{Finn_McOwen_IUMJ_93}.
The gap between such requirents and the 
sufficient condition of 
Loewner and Nirenberg
still remained broad.
In section~{1.3}
we show how
all these  results
follow from Theorem~\ref{maintheorem}.

This paper studies 
(\ref{singyam}),
(\ref{R=-1})
in dimensions
$n\geq 3$.
In the case of
$\mathbf{S}^2$
the complete  metric conformal to
$\gstand$
and having the constant negative curvature is
called the 
Poincare metric.
The equation for Poincare metric 
is slightly different from
(\ref{singyam}).
Mazzeo and Taylor
\cite{Mazzeo_Taylor_uniformisation} 
proved that 
the Poincare metric in
$\Omega\subset\mathbf{S}^2$
always exists provided the complement
$K$
has at least two distinct points.

In 1988 Schoen and Yau 
\cite{Schoen_Yau_Invent_88}
were led by their research 
on locally  conformally flat maniofolds
to the case
$R\geq 0$.
They found that a  necessary
condition for solvability of
(\ref{singyam})
with
$R\geq 0$,
as oposed
to 
the case
(\ref{R=-1}),
is smallness of
$K$.
For example, they proved that
if the solution exists then the Newtonian capacity of
$K$
must vanish.
They also established that the solvabilty  
in this case 
implies that
\begin{equation*}
\mathcal{H}^{\varepsilon+(n-2)/2}(K)=0\quad
{\rm
for
\quad 
all
\quad
}
\varepsilon>0.
\end{equation*}
Thus the Hausdorff dimension
$(n-2)/2$
separates
the cases of the negative and non-negative
curvature.
Similarly to the 
negative curvature case, the existence of
a solution is known at the moment only in cases when
$K$
has much more structure that vanishing Hausdorff measure
or capacity. Despite the similarity in statements,  
the  results in the case
$R=1$
are much more  difficult to prove.
In a seminal paper
\cite{Schoen_CPAM}
Schoen established the existence of 
(\ref{singyam})
with
$R=1$
when
$K$
is a finite number (at least two)  of points.
Mazzeo and Pacard
\cite{Mazzeo_Pacard_Duke}
generalising earlier results
\cite{Mazzeo_Smale_JDG},
\cite{Pacard_TMNA},
\cite{Mazzeo_Pacard_JDG}
extended Schoen's result to the case
when
$K$
is a finite number of disjoint smooth 
submanifolds of the dimension
$k\leq (n-2)/2$.
There is also a construction of a solution 
using Kleinian  groups
in
the case when 
$K$
is a certain Cantor-type set
\cite{Schoen_Yau_Invent_88}

The case
$R=0$
is easier
becuase equation 
(\ref{singyam})
becomes linear.
In this case
it is known that the solution exists
provided that 
$K$
is essentially a finite union
of Lipschitz submanifolds 
 of  dimension
$k\leq (n-2)/2$
\cite{Delanoe_Contemp_Math_90},
\cite{Ma_McOwen},
\cite{Kato_Nayatani}.

In the paper  we are interested only 
in the basic  problem  of existence for
(\ref{singyam}).
However,  other questions about solutions of
(\ref{singyam})
can be asked as well.
For example  problems of  
uniqueness, asymptotic behaviour of 
$u$ 
near
$\partial\Omega$,
structure of moduli space of solutions, gluing different solutions, 
are  investigated in
\cite{KMPS},
\cite{Mazzeo_Pollack_Uhlenbeck_JAMS},
\cite{Mazzeo_Pollack_Uhlenbeck_TPMNA},
\cite{McOwen_survey}.
Some of the results mentioned above hold for more general manifolds than 
$\mathbf{S}^n$.
The result directly related to the present paper was proved by Aviles 
and McOwen
\cite{Aviles_McOwen_Duke_88},
\cite{Aviles_McOwen_JDG_85},
\cite{Aviles_McOwen_JDG_88}.
They established that 
an open subset 
$\Omega$
of any closed Riemannian manifold
$(M,g)$
admits a complete metric with the constant negative scalar curvature
conformal to
$g$
provided
$K$
is a finite union of closed smooth submanifolds of dimensons
$k>(n-2)/2$.
In a future publication we 
introduce a suitable capacity and
extend
our Theorem~\ref{maintheorem}
to more general manifolds.


\subsection{Main theorem}
\label{sec1.2}

We investigate the solvability of
(\ref{singyam})
by attracting ideas from the nonlinear  potential theory.
Let us recall a fundamental result form the classical 
potential theory for the Laplace equation.
This is the Wiener test for  the classical 
Dirichlet problem for harmonic functions 
\cite{Wiener}.
Wiener theorem states that the Dirichlet 
problem 
$$
\left\{
\begin{array}{rcl}
\Delta w =0
&
{\rm in}
&
D
\\
w=f
&
{\rm on}
&
\partial D
\end{array}
\right.
$$
in a bounded domain 
$D\subset{\bf R}^n$,
$n\geq 3$,
is solvable for all
boundary data 
$f\in C(\partial D)$
if and only if
${\bf R}^n\setminus D$
is 
{\it not thin}. 
Explicitly the latter means that 
\begin{equation*}
\int_0^1
\frac{cap (B(x, r)\setminus D)}{cap (B(x, r))}
\,
\frac{dr}{r}
=+\infty
\quad
{\rm for \quad  any}
\quad
x\in\partial D.
\end{equation*}
Here
$1$
can be replaced by any
$\delta>0$,
and
$cap$
is the classical
(electrostatic) capacity.
Our main Theorem~\ref{maintheorem}
states that problem
(\ref{singyam}),
(\ref{R=-1})
admits a solution if and only if
a Wiener-type test with a 
certain capacity holds.

Let us scetch the definition of
the capacity apropriate for  problem
(\ref{singyam}), 
see 
section~\ref{capsection}
for more details.
Take a compact set 
$E\subset{\bf S}^n$,
$n\geq 3$,
with
\begin{equation*} 
\diam_{\gstand}(E) \leq \pi/3.
\end{equation*} 
After a rotation we can assume that 
such
$E$ 
lies in the southern hemisphere.
We set 
\begin{eqnarray*}
\capacity(E)
&
= 
&
\inf 
\left\{
\int_{{\bf S}^n}
\left|
\nabla^2 \varphi
\right|^{(n+2)/4}
\,
dvol_{\gstand}
\right\},
\\
&
&
\frac{1}{(n+2)/(n-2)}
+
\frac{1}{(n+2)/4}
=1.
\end{eqnarray*}
Here 
symbols
$dvol_{\gstand}$,
$\nabla$,
and
$| \cdot |$,
stand respectively
for the 
volume element,
connection,
and  
norm
with respect to 
the metric
$\gstand$.
The infimum  is taken over all
$\varphi\in C^\infty({\bf S}^n)$
such that
$\varphi|_{E} \geq 1$
and
$\varphi 0$
on
the 
northern hemisphere.
Essentially,
$\capacity$
is the  
{\it Bessel capacity}
for the Sobolev space
$W^{2, (n+2)/4}({\bf R}^n)$.
Bessel capacities have been intensively
investigated in the nonlinear potential theory.
Nonlinear potential theory originates in early works of 
Maz'ya and Serrin in the 1960s
and was extensively  developed later in 1970s and 1980s
by many authors.
Our paper heavily relies on it.
The main references will be monographs
by Adams and Hedberg
\cite{Adams_Hedberg_book},
Maz'ya
\cite{Maz'ya_book},
and
Ziemer
\cite{Ziemer_book}.
There  the reader can also find
a rich bibliography and historical notes
Now we state the main theorem.

\begin{theorem}
\label{maintheorem}
Let
$\Omega\subset{\bf S}^n$,
$n\geq 3$,
be an open set and
$K={\bf S}^n \setminus \Omega$.
Then the following properties are equivalent:

\noindent
(i)
In 
$\Omega$
there exists a complete metric 
with constant negative scalar curvature
conformal to 
$\gstand$.

\noindent
(ii)
The compactum
$K$
is not thin, that is
for  any
$p\in K$ 
\begin{equation}
\label{wienertest}
\int_0^{1/2}
\left(
\frac{\capacity(B(p, r)\cap K)}{\capacity(B(p,r))}
\right)^{2/(n-2)}
\,
\frac{dr}{r}
=+\infty.
\end{equation}
\end{theorem}

Wiener test
(\ref{wienertest})
is a capacitary condition on
$K$.
Geometric  properties of
the capacity
$\capacity$
are well understood
due to  investigations
in nonlinear potential theory.
Using the information avialable there,
we show in 
section~\ref{sec1.3}
that more transparent geometric results 
easily follow from
Theorem~\ref{maintheorem}.

In view of 
Theorem~\ref{maintheorem}
it would be  interesting to clarify
how   the condition
\begin{equation*}
\capacity({\bf S}^n\setminus \Omega)=0
\end{equation*}
relates to the 
conformal deformation to 
nonnegative scalar curvature
$R\geq0$.
In 
\cite{Labflat}
we apply potential theory ideas to the scalar flat 
case
$R=0$. 
In this situation 
as opposed to 
Theorem~\ref{maintheorem}
the set
$K$
should be small.

We mention that ideas from potential theory 
have been  used in conformal geometry before.
For example
Schoen and Yau
\cite{Schoen_Yau_Invent_88},
\cite{Schoen_Yau_book}
used capacity related to Sobolev space
$W^{1, q}({\bf R}^n)$.
Capacity
$\capacity$
was implicitely used 
at some stage
in
\cite{KMPS}
to prove distribution removability 
of isolated singularities
for  the equation in
(\ref{singyam})
with
$R=1$.

\begin{remark}
\label{rem}
Intuitevely,
the completeness condition
forces solutions of
(\ref{singyam})
to blow up in some sense near
$\partial\Omega$.
Dhersin and LeGall
\cite{DLeG}
considered the problem
\begin{equation}
\label{largeproblem}
\left\{
\begin{array}{rcl}
\Delta u - u^2 =0
&
{\rm in}
&
D
\\
u(x)\to+\infty
&
{\rm when}
&
x\to\partial D
\end{array}
\right.
\end{equation}
in general domains
$ D \subset{\bf R}^n$.
They proved
that a Wiener test 
characterises domains
$D$
for which
(\ref{largeproblem})
is solvable.
Their paper provided an important inspiration for 
our work, although
the consideration in 
\cite{DLeG}
is based on 
probabilistic methods.
In fact,
there is a strong  connection between 
$u$
in
(\ref{largeproblem})
and 
a certain branching random process
(so-called  {\it Brownian snake})
\cite{LeGall_book},
\cite{Dynkin_book}.
To find an adequate probabilistic interpretation
for
$p>2$
is an important open problem in the area
\cite{LeGcongress},
\cite{LeGall_book}
\cite{Dynkin_book}.
However, in
\cite{Lab2}
we established a 
solvability criterion
for
\begin{equation*}
\Delta u - u^p =0
\end{equation*}
blowing up at the boundary
in the full range of
$p>1$.
Relying there on entirely analytic ideas
we establish 
estimates for 
solutions 
in terms of  the capacity
associated with the variational integral
\begin{equation*}
\int_{{\bf R}^n}|D^2\varphi|^{p'},
\quad
\varphi\in C^\infty_0({\bf R}^n),
\quad
\frac{1}{p}+\frac{1}{p'}=1.
\end{equation*}
Comparison of 
Theorem~\ref{maintheorem}
with the condition from
\cite{Lab2}
implies that
any conformal factor in 
(\ref{singyam}),
(\ref{R=-1})
must blow up pointwisely near
$K$.
Iterplay between Brownian motion on a Riemannian manifold
and geometric properties of the manifold is an important subject.
We refer to Grigor'yan's survey 
\cite{Grigor'yan}
for an exposition of the classical and new results.
It would be interesting to understand the relation between
super-Brownian motion, Brownian snakes
and the geometric problems from  this paper (or properties of Riemannian 
manifolds in general).
Some results in this direction can be found in
\cite{Dynkin_process_manifolds}.
\end{remark}

\subsection{Examples}
\label{sec1.3}

We illustrate how Theorem~\ref{maintheorem}
allows  to establish the existence for
(\ref{singyam}),
(\ref{R=-1})
in concrete situations. In particular,
apparently
all necessary or sufficient conditions 
from previous papers can be easily derived from
(\ref{wienertest}).
The reason for this
is that 
capacity
$\capacity$
had 
appeared before
in different problems related to the interaction between 
nonlinear potentials
and the Littlewood-Paley
theory.
As a result, it was intensively studied in the 1970s-1980s and
its geometric
properties 
are well known.
They are carefully documented e.g. in
\cite{Adams_Hedberg_book},
\cite{Maz'ya_book}
\cite{Ziemer_book}.

\begin{example}
The necessity of Loewner-Nirenberg condition
\cite{Loewner_Nirenberg}
\begin{equation*}
\Hmeas^{(n-2)/2}(K)=+\infty
\end{equation*}
for solvability of
(\ref{singyam}),
(\ref{R=-1})
follows immediately
from 
Theorem~\ref{maintheorem}
and the implication
\begin{equation*}
\Hmeas^{(n-2)/2}(K)<+\infty
\Longrightarrow\capacity(K)=0,
\end{equation*}
valid for our capacity
$\capacity$,
see 
\cite{Adams_Hedberg_book} and
(\ref{capHausdmeas})
in
section~\ref{capsection}
below.
\end{example}

\begin{example}
Assume that
$K$
is a smooth {\it immersed}  submanifold 
in
${\bf S}^n$
of dimension
$k$,
\begin{equation}
\label{nice_dimension}
k>(n-2)/2.
\end{equation}
Then
(\ref{singyam}),
(\ref{R=-1})
has a solution.
In fact,  
according  to 
Theorem~\ref{maintheorem}
we need to show that
(\ref{wienertest})
holds.
Fix
any
$p\in K$.
By the definition of immersion
there exists an open  
smooth submanifold
$E$
of dimension
$k$
embedded in
${\bf S}^n$
such that
$p\in E$
and
$E\subset K$.
Therefore
\begin{equation*}
\capacity(K\cap B(p,r))
\geq
\capacity(E \cap B(p,r)).
\end{equation*}
The exponential map
$\exp_p$
is a diffeomorphism of a neighbourhood
of the origin
$0$ 
in
$T_p\mathbf{S}^n$.
In the sufficiently small neighbourhood of
$p$
the smooth submanifold
$E$
is well approximated by 
the image under
$\exp_p$
of 
a 
neighbourhood
of
$0$
in
$T_pE$,
$T_pE\subset T_p\mathbf{S}^n$.
Capacities of a set and its image 
under a diffeomorphism are equivalent
\cite{Adams_Hedberg_book}.
Hence
utilising the 
scaling property 
(\ref{capscaling})
we find a small number
$r_0>0$,
such that
\begin{equation*}
\capacity(E \cap B(p,r))
\geq C(E)
\left(
\frac{r}{r_0}
\right)
^{(n-2)/2} 
\capacity(E\cap B(p, r_0))
\quad
{\rm for}
\quad
r\in (0, r_0).
\end{equation*}
Here the constant
$C(E)>0$
depends on the smoothness of 
$E$.
Capacity of the ball can
be estimated by  
\begin{equation*}
\capacity(B(p, r))
\asymp
r^{(n-2)/2}
\quad
{\rm for}
\quad
r \in (0, r_0),
\end{equation*}
see  
(\ref{capball}).
Now
(\ref{wienertest})
follows.
\end{example}

\begin{example}
Let
$d=\varepsilon + (n-2)/2$,
$\varepsilon>0$.
Assume that for any
$x\in K$
there is a positive constant
$C$
such that
\begin{equation}
\label{density_cond}
\mathcal{H}^d_\infty(K\cap B(x,r))\geq Cr^d
\end{equation}
for all
$r$
near 
$0$.
That is, the Hausdorff $d$-content of
$K$
is big
at all small scales.
Then
the conformal metric from
Theorem\ref{maintheorem}
exists. This follows at once from
(\ref{capHausdcontent}),
(\ref{capball}).
Density condition 
(\ref{density_cond})
allows to recover 
the results of Finn
\cite{Finn_94},
\cite{Finn_99},
\cite{Finn_00}
about sets
$K$
with stratified cone-type
tangent structure.
Also 
(\ref{density_cond})
allows to establish
existence of the complete metrics
in the cases when
$K$
satisfies different 
geometric conditions
invariant with respect to 
quasiconformal maps,
see e.g.
\cite{Koskela}.
\end{example}

\begin{example}
Let 
$K$ 
be the Lebesgue cusp.
That is for a fixed
$\rho>0$
and for a
continuous positive nondecreasing function
$h$
on the real line,
$h(r)=O(r)$,
$r\to 0$,
we set
\begin{equation*}
T_h
=
\left\{
x\in\mathbf{R}^n
\colon
0\leq x_n\leq \rho,
\
\left(
x_1^2+\cdots+x_{n-1}^2
\right)^{1/2}
\leq h(x_n)
\right\}.
\end{equation*}
Then define
$K$
to be the preimage of
$T_h$
under the stereogrphic projection,
\begin{equation*}
K=\sigma^{-1}(T_h).
\end{equation*}
The existence
of the singular conformal metric
$g$
from 
Theorem~\ref{maintheorem}
in
$\mathbf{S}^n \setminus K$
depends on the dimension
$n$.
If 
$n=3$
then
$g$
always exists.
For higher dimensions
$g$
exists if and only if
\begin{eqnarray*}
&
&
\int_0^1
\frac
{dr
}
{
r\log
\left(r/h(r) \right)
}
=+\infty
\quad
{\rm for}
\quad
n=4
\\
&
&
\int_0^1
\left(
\frac{h(r)}{r}
\right)^{(n-4)/(n-2)}
\,
\frac{dr}{r}
=+\infty
\quad
{\rm for}
\quad
n>4.
\end{eqnarray*}
Indeed, we only need to check
that
(\ref{wienertest})
holds for the South pole
$S=\sigma^{-1}(0)$.
To  verify 
(\ref{wienertest}),
first recall
that for the cyllinder
\begin{equation*}
\Pi=\sigma^{-1}
\Big(
(-\delta, \delta)
\times 
\cdots 
\times 
(-\delta, \delta)
\times
(r/2, r)
\Big),
\quad
4\delta<r,
\end{equation*}
with
$r>0$
small enough,
the capacity 
is given by the following
formulae
\cite{Maz'ya_book}, Ch. 9:
\begin{eqnarray*}
\capacity( \Pi )
&
\asymp
&
\capacity( B(S, r))
\quad
{\rm for}
\quad
n=3,
\\
\capacity( \Pi )
&
\asymp
&
\frac
{1}
{\log\left( r/\delta \right)}
\capacity( B(S, r))
\quad
{\rm for}
\quad
n=4,
\\
\capacity( \Pi )
&
\asymp
&
\left(
\frac{\delta}{r}
\right)^{(n-4)/2}
\capacity( B(S, r))
\quad
{\rm for}
\quad
n>4.
\end{eqnarray*}
Now just  apply elementary estimate
(\ref{sertointegr}).
\end{example}

\subsection{Organisation of the paper}
\label{sec1.4}

In section 2 we 
introduce the capacity  and use it to  
prove some preliminary estimates for solutions
of the equation.
We also describe there the unique feature of 
equation
(\ref{singyam}),
(\ref{R=-1}).
Namely,
the existence of a  finite maximal solution
$u_\Omega$
dominating all other solutions pointwisely.

In section 3
we prove the crucial estimates for solutions of
(\ref{singyam}),
(\ref{R=-1})
in terms of the capacity.
The principal  difficulty in the proof of the main
Theorem~\ref{maintheorem}
is the analysis of the completeness condition in
(\ref{singyam}).
This involves understanding the behaviour 
of the conformal factor
$u$
near
$\partial\Omega$
under no assumptions
(say, when proving
(i)$\Rightarrow$(ii)
in Theorem~\ref{maintheorem})
on the structure of
$\Omega$.
The first estimate in section 3,
Theorem~\ref{captheorem},
controlls the solution pointwisely 
away from
$\partial\Omega$.
The second estimate, Theorem~\ref{intu2/(n-2)},
provides the integral control
when we stay arbitrarily close to
$\partial\Omega$.

With all this background in place
we proceed to prove our main 
result,
Theorem\ref{maintheorem},
in section 4.
The sufficiency of the Wiener test
(\ref{wienertest})
will follow rather straighforwardly from
the pointwise 
estimate from section 3.
To prove the necessity we supose 
that the negation  of
(\ref{wienertest})
holds. In other words,
suppose that 
the complement of
$\Omega$
is thin at some point.
We will find  a curve 
in
$\Omega$
approaching this point, 
such that its length with respect to
$u_\Omega$
(and hence with respect to any other solution of
(\ref{singyam}),
(\ref{R=-1}))
is finite.
How to construct such a curve without any assumptions on
$\partial\Omega$?
The key idea here is to reduce this issue to an integral estimate.
To achieve this we bring in the estimate from Theorem~\ref{intu2/(n-2)}.

Throughout this paper,
we will use the notation
\begin{equation*}
q=\frac{n+2}{n-2},
\quad
q'=\frac{n+2}{4},
\quad
\frac{1}{q}+\frac{1}{q'}=1.
\end{equation*}
By
$\geucl$
we
denote the Euclidean metric in
$\mathbf{R}^n$ and by
$\gstand$
the standard meric on the sphere
induced by
$\geucl$.
By
$B(p,r)$
we denote the ball
of radius
$r$
centered at 
$p$
for
$\gstand$
or
$\geucl$. 
It will be clear from the context 
which metric is taken.
For an integer
$j$
we put
$r_j=2^{-j}$.
By
$B_j$
we denote the dyadic
ball
in
$\mathbf{R}^n$,
$B_j=B(0,r_j)$.
We denote
the Green's function for the Laplacian
in
$B(0,R)\subset\mathbf{R}^n$
by
$G_R$.
By
$C$,
$\widetilde{C}$,
$C_1$,
$\dots$,
we denote positive constants depending only on the dimension.
The value of
$C$,
$\widetilde{C}$,
$C_1$,
$\dots$,
may vary even within the same line.
We write
\begin{equation*}
A\lesssim B
\
(A\gtrsim B)
\end{equation*}
if
\begin{equation*}
A\leq CB
\
(A\geq CB)
\end{equation*}
for some
$C$.
We write
\begin{equation*}
A\asymp B
\end{equation*}
if
\begin{equation*}
A\lesssim B\lesssim A.
\end{equation*}

\subsection{Acknowledgements.}

I would like to thank Rick Schoen
for many  intersting and stimulating 
discussions on the problem.
I also would like to thank
Michael Struwe, Tom Ilmanen, Reiner Sch\"atzle,
Neil Trudinger, and Xu-Jia Wang
for their interest and support.

\section{Preliminaries on the  equation and capacity}
\label{sec2}

\setcounter{equation}{0}


\subsection{Equation on sphere and in space}
\label{sphere&Rn}

Let
$g$
be a metric on a manifold 
$M$
of dimension
$n$,
$n\geq 3$.
The operator
\begin{equation*}
\mathcal{L}_g
=
-
4\frac{n-1}{n-2}
\Delta_g 
+
R(g)
\end{equation*}
from
(\ref{scal_curv_change})
is called conformal Laplacian
\cite{Schoen_Yau_book}.
If we change the metric conformally
\begin{equation*}
\hat{g}=\varphi^{4/(n-2)}g,
\end{equation*}
then
\begin{eqnarray*}
&
&
R(\hat{g}) = \varphi^{-(n+2)/(n-2)} 
\mathcal{L}_g \varphi,
\\
&
&
\mathcal{L}_{\hat{g}} v 
= 
\varphi^{-(n+2)/(n-2)} \mathcal{L}_g \left(\varphi \, v \right).
\end{eqnarray*}
More generally, 
let 
$\tilde{M}$
be  another manifold with
the metric
$\tilde{g}$, 
and let
$f\colon M\to \tilde{M}$
be a diffeomorphism.
Assume that 
$f$
changes the metric conformally
\begin{equation*}
f^*\tilde{g}=\varphi^{4/(n-2)} g.
\end{equation*}
Then 
\begin{eqnarray*}
&
&
f^*(R(\tilde{g})) = \varphi^{-(n+2)/(n-2)} 
\mathcal{L}_g \varphi,
\\
&
&
f^*(\mathcal{L}_{\tilde{g}} v) = 
\varphi^{-(n+2)/(n-2)} \mathcal{L}_g \left(\varphi \, f^*v\right).
\end{eqnarray*}
Now, 
the stereographic projection
$
\sigma\colon \mathbf{S}^n \setminus \{N\}
\to \mathbf{R}^n
$
is a conformal diffeomorphism between
$(\mathbf{S}^n \setminus\{N\}, \gstand)$
and
$(\mathbf{R}^n, \geucl)$
because
\begin{eqnarray*}
\left(
\sigma^{-1}
\right)^*
\gstand
&
=
&
\left(
\frac{2}{1+|x|^2}
\right)^2
\geucl
\\
&
=
&
\Upsilon^{4/(n-2)}
\geucl
\end{eqnarray*}
with
\begin{equation*}
\Upsilon(x)
=
\left(
\frac{2}{1+|x|^2}
\right)^{(n-2)/2}
\quad
x\in\mathbf{R}^n.
\end{equation*}
According to the above formulae for
the conformal changes
we have the following correspondence.

Let
$\Omega\in\mathbf{S}^n$,
$N\subset\Omega$,
and let the function
$v$
satisfy
\begin{equation}
\label{eq_on_sphere}
v^{-(n+2)/(n-2)}
\mathcal{L}_{\gstand}
v 
=
-
1,
\quad
v>0
\quad
{\rm in}
\quad
\Omega.
\end{equation}
Then
the function
\begin{eqnarray}
\label{conf_change_of_conf_factor}
u(x)
&
=
&
\Upsilon(x)
\left(
\sigma^{-1}
\right)^*v (x)
\nonumber
\\
&
=
&
\Upsilon(x)
v(\sigma^{-1}x),
\quad
x\in\mathbf{R}^n,
\end{eqnarray}
satisfies
\begin{equation*}
u^{-(n+2)/(n-2)}
\mathcal{L}_{\geucl} u 
=
-
1,
\quad
u>0
\quad
{\rm in}
\quad
\sigma(\Omega).
\end{equation*}
Thus after multiplication
by a constant, 
$u$
satisfies
\begin{equation}
\label{positsoltomaineq}
u>0,
\quad
\Delta u 
- 
u^
{
(n+2)/(n-2)
} 
=0.
\end{equation}
Conversely,
for any  solution
$u$
of
(\ref{positsoltomaineq})
defined
in
$\sigma(\Omega)$
set
\begin{equation*}
v=\sigma^* 
\left(
u
/\Upsilon
\right)
.
\end{equation*}
Then after  multiplication  by a constant,
$v$ 
satisfies
(\ref{eq_on_sphere})
in
$\Omega\setminus \{ N \}$.
Moreover 
we know 
\cite{Loewner_Nirenberg}
that
for any
$u$
solving
(\ref{positsoltomaineq})
in a neighbourhoud of infinity
in
$\mathbf{R}^n$,
there exists a constant
$A>0$
such that
\begin{equation*}
u(x)
=
\frac{A}{|x|^{n-2}}
+
o
\left(
\frac{1}{|x|^{n-2}}
\right),
\quad
x\to\infty.
\end{equation*}
Hence
we can
extend
$v$
to
$N$
by continuity, 
remove the isolated singularity,
and conclude that
(\ref{eq_on_sphere})
holds.

Clearly
the metric
$v^{4/(n-2)}\gstand$
is complete in
$\Omega$
for
$v$
from 
(\ref{eq_on_sphere})
if and only if
$u^{4/(n-2)}\geucl$
is complete in
$\sigma(\Omega)\cup\{\infty\}$
for the corresponding
$u$
from
(\ref{conf_change_of_conf_factor}),
(\ref{positsoltomaineq}).

The previous discussion shows that
the existence of the singular Yamabe metric in a domain
on the sphere is equivalent to finding a complete solution of
(\ref{positsoltomaineq})
in the exterior domain
in
$\mathbf{R}^n$.
Let us describe the main features
of equation
(\ref{positsoltomaineq})
in
$\mathbf{R}^{n}$.
Omited proofs  can be found for example in
\cite{Loewner_Nirenberg}.

The crucial fact  about  solutions of
(\ref{positsoltomaineq})
that  will be used constantly in this  paper 
is the 
{\it elliptic comparison principle}.
As a consequence of this principle, 
local regularity estimates 
hold for 
$u$.
In particular, if
$u\in L^q_\loc$
is a distributional  solution of
(\ref{positsoltomaineq})
then,
in fact,
$u\in C^\infty_\loc$
and
$u$
is the classical solution.
Moreover, let
$u$
be 
{\it any} 
solution of
(\ref{positsoltomaineq})
in
an open set
$O\subset\mathbf{R}^n$.
Then
\begin{equation}
\label{KellOsserm}
u(x)
\lesssim
\frac{1}
{
\dist(x,\partial O)
^{(n-2)/2}
}
\quad
{
\rm 
for
\quad
all 
\quad
}
x\in O.
\end{equation}
This is an estimate uniform in
$u$. 
It was first discovered
by Keller
\cite{Kell}
and 
Osserman
\cite{Osserm},
and also follows from the comparison principle.

Estimate
(\ref{KellOsserm})
combined with the elliptic
Perron argument implies
the existence of the finite
solution
$u_O$
which is {\it maximal}
in
$O$.
It means that
the inequality
\begin{equation*}
u\leq u_O
\quad
{\rm in}
\quad
O
\end{equation*}
holds  for any   other 
$u$
solving
(\ref{positsoltomaineq})
in
$O$.
Clearly
\begin{equation*}
u_{O_1}\leq u_{O_2}
\quad
{\rm in}
\quad
O_2
\quad
{\rm when}
\quad
O_1 \supset O_2.
\end{equation*}
Let
$K_1$, 
$\ldots$, 
$K_m$
be compact sets
in
$\mathbf{R}^n$,
let
\begin{equation*}
K=K_1\cup\cdots\cup K_m,
\end{equation*}
and let
$u$, 
$u_1$, 
$\ldots$, 
$u_m$
be the maximal solutions of
(\ref{positsoltomaineq})
in
$K^c$, 
$K_1^c$, 
$\ldots$, 
$K_m^c$
respectively.
Then the 
H\"older inequality
and the comparison ensure that
\begin{equation}
\label{sum}
m^
{
-(n-2)/(n+2)
}
\sum_{i=1}^m u_i
\leq
u
\leq
\sum_{i=1}^m u_i
\quad
{\rm in}
\quad
K^c.
\end{equation}
If
$x_0\in\partial O$
and
$(\partial\Omega)\cap B(x_0, r)$
is a smooth hupersurface for some
$r>0$,
then
\begin{eqnarray}
\label{asympformaxsol}
u_O(x)
\,
\dist(x,\partial  O)^{(n-2)/2}
&
\to
& 
\left(
\frac
{
n(n-2)
}
{
4
}
\right)^{(n-2)/4}
,
\nonumber
\\
&
&
{\rm when}
\quad
x\to(\partial\Omega)\cap B(x_0, r).
\end{eqnarray}
Asymptotic behaviour
(\ref{asympformaxsol})
holds  in fact for 
{\it any}
solution of
(\ref{positsoltomaineq}) 
blowing up at
$(\partial\Omega)\cap B(x_0, r)$.

Solutions to
(\ref{positsoltomaineq})
exhibit the following dilation invariance:
for all
$a>0$
and 
$r>0$,
\begin{eqnarray}
\label{scalingeq}
u
\
{\rm solves}
\
(\ref{positsoltomaineq})
\
{\rm in}
\
B(0,r)
&
\Longrightarrow
&
a^{(n-2)/2}u(a\cdot)
\
{\rm solves}
\
(\ref{positsoltomaineq})
\nonumber
\\
&
&
{\rm in}
\
B(0, r/a).
\end{eqnarray}

Finally consider equation
(\ref{eq_on_sphere})
on the sphere.
As a direct consequence of the properties of the 
stereographic projection
which we discussed above,
there exists the maximal solution of
(\ref{eq_on_sphere})
and the estimates analogous to
(\ref{KellOsserm})--(\ref{asympformaxsol})
hold.

\subsection{Capacity}
\label{capsection}

In this paragraph
we define the capacity 
$\capacity$
for subsets of the unit sphere
$\mathbf{S}^n$.
Esentially it is a particular Bessel capacity
$\capinRn$
in
$\mathbf{R}^n$.
Omited proofs of the statements
about
$\capinRn$
can be found in monographs
\cite{Adams_Hedberg_book},
\cite{Maz'ya_book},
and
\cite{Ziemer_book}.
By
$\mathbf{S}^n_S$
and
$\mathbf{S}^n_N$
we denote
southern and northern hemispheres.

Fix the spherical cup
around the south pole
$S$
by writing
\begin{equation}
\label{spherical_cup}
\mathbf{U}
=
\left\{
x\in\mathbf{S}^n\colon
d_{\gstand}
(S,x)\leq \pi/3
\right\},
\quad
\mathbf{U} \subset \mathbf{S}^n_S.
\end{equation}
Take any compact set
$K\subset \mathbf{S}^{n}$
with
\begin{equation*}
\diam_{\gstand} (K)\leq\pi/3
.
\end{equation*}
The rotation  group
$SO(n+1)$
acts transitively on
$\mathbf{S}^{n}\hookrightarrow\mathbf{R}^{n+1}$.
Map 
$K$ 
by a rotation
$\Phi\in  SO(n+1)$
in a way that
$\Phi(K)\subset\mathbf{U}$.
Define
\begin{eqnarray}
\label{capinsphere1}
\capacity(K)
=
\inf
\Big\{
\int_{\mathbf{S}^{n}}
\left|
\nabla^{2}
\varphi
\right|^{(n+2)/4}
\,
dvol_{\gstand}
\colon
&
&
\varphi\in C^{\infty}(\mathbf{S}^{n}),
\
\varphi |_{\mathbf{S}^{n}_{N}}=0,
\nonumber
\\
&
&
\varphi |_{\Phi(K)}\geq 1
\Big\}.
\end{eqnarray}
We will prove  that different choices of
$\Phi\in SO(n+1)$
lead to equivalent
capacities. First we give an alternative
description of the capacity.
Stereographic projection
$\sigma$  
is a smooth  quasiisometry
between
$\mathbf{S}^{n}_{S}$
and
$B(0,1)\subset\mathbf{R}^n$.
Hence
\begin{eqnarray}
\label{capinsphere2}
\capacity(K)
\asymp
\inf
\Big\{
\int_{\mathbf{R}^{n}}
\left|
D^{2}
\psi
\right|^{(n+2)/4}
\,
dx
\colon
&
&
\psi\in C^{\infty}_{0} 
(B(0,1)),
\nonumber
\\
&
&
\psi |_{\sigma\circ\Phi(K)}\geq 1
\Big\}.
\end{eqnarray}
Let us introduce the corresponding capacity for sets in
$\mathbf{R}^{n}$. 
For a compact set
$E\subset B((0,1)) \subset\mathbf{R}^{n}$
its Bessel capacity is defined as
\begin{eqnarray}
\label{capinRn}
\capinRn(E)
=
\inf
\Big\{
\int_{\mathbf{R}^{n}}
\left|
D^{2}
\psi
\right|^{(n+2)/4}
\,
dx
\colon
&
&
\psi\in C^{\infty}_{0} 
(B(0,2)),
\nonumber
\\
&
&
\psi |_E 
\geq 1
\Big\}.
\end{eqnarray}
Notice that the set
$\sigma\circ\Phi(K)$
stays away from the boundary of the unit ball:
\begin{equation*}
\sigma\circ\Phi(K) 
\subset
B(0, 99/100).
\end{equation*}
Hence properties of Bessel capacities imply  
that 
for
$E=\sigma\circ\Phi(K)$
the right hand sides of
(\ref{capinRn})
and
(\ref{capinsphere2})
are equivalent.
Thus
\begin{equation}
\label{capequivalence}
\capacity(K)
\asymp
\capinRn (\sigma\circ\Phi(K) ).
\end{equation}
Now take another 
${\tilde{\Phi}\in SO(n+1)}$,
$\tilde{\Phi}(K)\subset\mathbf{U}$.
The same variational procedure as
(\ref{capinsphere1})
gives
the new capacity
$\widetilde{\capacity}(K)$.
Bessel capacity
(\ref{capinRn})
of a compactum and of its image under a bi-Lipschitz
homeomorphism are equivalent.
Apply  this to the 
locally bi-Lipschitz map
\begin{equation*}
\sigma
\circ
\tilde{\Phi}
\circ
\Phi^{-1}
\circ
\sigma^{-1}
\colon
\mathbf{R}^n
\to
\mathbf{R}^n
\end{equation*}
which sends
$\sigma\circ\Phi(K)$
to
$\sigma\circ\tilde{\Phi}(K)$, 
and utilise
(\ref{capequivalence})
to derive that
\begin{equation*}
\capacity(K)
\asymp
\widetilde{\capacity}(K).
\end{equation*}

Clearly property 
(\ref{wienertest})
of the set to be 
not thin 
does  not change when we pass to an
equivalent capacity.
Set functions
$\capacity$
and
$\capinRn$ 
enjoy subadditvity
and monotonicity
properties.
Standard scheme of axiomatic potential theory
extends them
to arbitrary sets
as the  outer measure.

Now we list some well-known metric estimates for
the capacity.
The following important scaling holds:
\begin{equation}
\label{capscaling}
\capinRn(tE)
\asymp
t^{(n-2)/2}
\capinRn(E),
\quad
t\in(0,1),
\quad
E\subset\subset B(0,1).
\end{equation}
Next,
for
$\alpha>0$
the 
{\it Hausdorff
$\alpha$-content}
of 
$E\subset\mathbf{R}^{n}$
(or $E\subset\mathbf{S}^{n}$)
is defined as
\begin{equation*}
\Hmeas^{\alpha}_{\infty}
(E)
=
\inf
\sum_{j}
r_j^{\alpha},
\end{equation*}
where the infimum is taken over all coverings
of
$E$
by countable unions of euclidean
balls
$\{B(x_{j}, r_{j})\}$
in
$\mathbf{R}^{n}$
(or
$d_{\gstand}$-balls
in
$\mathbf{S}^{n}$).
The set function
$\Hmeas^{\alpha}_{\infty}$
is subadditive and monotone.
For the Hausdorff measure we have
\begin{equation*}
\Hmeas^\alpha(E)=0
\Longleftrightarrow
\Hmeas^\alpha_\infty(E)=0.
\end{equation*}
There is a strong connection between
the capacity and the Hausdorf content and measure.
For any 
\begin{equation*}
\alpha >(n-2)/2
\end{equation*}
there is a constant
$C(n,\alpha)>0$
such that
\begin{equation}
\label{capHausdcontent}
\left(
\Hmeas^{\alpha}_{\infty}
(E)
\right)^{(n-2)/2}
\leq
C(n,\alpha)
\capinRn(E)^{\alpha},
\quad
E\subset\subset B(0,1).
\end{equation}
Hence sets of the capacity 
$0$
have  the Hausdorf dimension at most
$(n-2)/2$.
In the converse direction
the following
implication holds
for the 
Hausdorf measure:
\begin{equation}
\label{capHausdmeas}
\Hmeas^{(n-2)/2}
(E)<+\infty
\Longrightarrow
\capinRn(E)=0
\end{equation}
for 
$E\subset\subset B(0,1)$.
According to  
(\ref{capequivalence})
statements
(\ref{capHausdcontent}) 
and
(\ref{capHausdmeas})
also hold for
$\capacity$.
From
(\ref{capequivalence}),
(\ref{capscaling}),
and
(\ref{capHausdcontent})
we also derive that
\begin{equation}
\label{capball}
\capacity(B(p,r))
\asymp
r^{(n-2)/2},
\quad
p\in\mathbf{S}^{n},
\quad
0\leq r\leq \pi/6.
\end{equation}

\subsection{An estimate}

In this paragraph we
provide an
integral  estimate for any  solution
$u$
of 
(\ref{positsoltomaineq})
outside  a compact
set
$K$
in terms of the capacity of 
$K$.
It will be frequently used in the sequel.
More precisely,
the following lemma produces a cut-off
function
$\eta$
which
vanishes in a neighbourhoud
of 
$K$,
equals
$1$
away from 
$K$,
and  bounds
the rate of a possible blow-up
of
$u$
via estimates
(\ref{eta})
and
(\ref{upowerqtimeseta}).

\begin{lemma}
\label{mainlemma}
Let
$K\subset B(0, 1)$
be a compact set
in
$\mathbf{R}^n$,
$n\geq 3$,
and 
\begin{equation*} 
m\geq\frac{n+2}{2}.
\end{equation*}
Let 
$u$
solve
(\ref{positsoltomaineq})
in
$K^c$.
Then there exists a  function
$\varphi\in C^\infty_0 (B(0,2))$
such that
$0\leq \varphi \leq 1$
in
$B(0,2)$,
$\varphi =1$
in an open neighbourhood of 
$K$,
\begin{equation}
\label{phi}
\int_{B(0,2)} |D^2\varphi|^{(n+2)/4}
\lesssim
\capinRn (K), 
\end{equation}
and such that for
$\eta=(1-\varphi)^m$
the inequalities
\begin{equation}
\label{eta}
\int_{{\bf R}^n}
u(|D\eta| +|\Delta\eta|)
\leq
C(m,n)
\capinRn (K),
\end{equation}
\begin{equation}
\label{upowerqtimeseta}
\int_{{\bf R}^n}
u^{(n+2)/(n-2)}\eta
\leq
C(m,n)
\capinRn (K)
\end{equation}
hold.
\end{lemma}

\begin{proof}
{\bf 1.}
The open set
$K^c$
can  be approximated from the interior
by domains  with smooth boundaries.
Consequently, by standard
continuity properties of  capacity,
we can assume 
in the proof
that
$K$
is a disjoint union of 
a finite number of closed domains with 
smooth boundaries.
We set
$\widehat{B}=B(0,2)$.

We claim that
there exists
a function
$\varphi\in C^\infty_0 (\widehat{B})$
with
$0\leq \varphi \leq 1$
in
$\widehat{B}$
and
$\varphi =1$
in an open neighbourhood of 
$K$
such that
(\ref{phi})
holds. 
To prove this, we first recall
a well-known result in 
nonlinear potential theory
\cite{Adams_Hedberg_book}
Chapter 2,
\cite{Maz'ya_book}
Chapter 9,
that states that
there exists a function
$\widetilde{\varphi}\in C^\infty_0(\widehat{B})$
such that
\begin{equation*}
\widetilde{\varphi}|_K\geq 1,
\quad
\int_{\widehat{B}}
|D^2\widetilde{\varphi}|^{q'}
\asymp \capinRn (K),
\quad
{\rm and}
\quad
\|\widetilde{\varphi}\|_{L^\infty(\widehat{B})}\lesssim 1.
\end{equation*}
Next, take  a function
$H\in C^\infty({\bf R}^1)$ 
such that
$$
H(t)=0
\quad
{\rm for}
\quad
t<1/3,
\quad
H(t)=1
\quad
{\rm for}
\quad
t>1/2.
$$
Now we take
$\varphi$
to be the smooth truncation
of
$\widetilde{\varphi}$,
$
\varphi
=
H(\widetilde{\varphi})
$.
Then
$$
\int_{\widehat{B}}
|D^2 \varphi|^{q'}
\lesssim
\int_{\widehat{B}}
|H''(\widetilde{\varphi})|^{q'}
|D\widetilde{\varphi}|^{2q'}
+
\int_{\widehat{B}}
|H'(\widetilde{\varphi})|^{q'}
|D^2 \widetilde{\varphi}|^{q'}
.
$$
To obtain
(\ref{phi}),
we just  apply the 
Gagliardo-Nirenberg interpolation inequality
\cite{Maz'ya_book}, Chapter 9,
to the first term:
if 
$1<r<\infty$,
then
for any
$f\in C_0^\infty (\widehat{B})$
\begin{equation}
\label{interpolation}
\| D f \|_{L^{2r}(\widehat{B})}
\lesssim
\| D^2 f \|_{L^{r}(\widehat{B})}^{1/2}
\| f \|_{L^\infty (\widehat{B})}^{1/2}
.
\end{equation}
We remark that   arguments 
of this type
are well known, cf.
\cite{Maz'ya_book}
Chapter 9,
\cite{Adams_Hedberg_book}
Chapter 3.

{\bf 2.}
Let
$u$
be a solution of
(\ref{positsoltomaineq}).
Take any
$\varepsilon>0$.
Appealing to decay 
(\ref{KellOsserm}),
we  choose
$R=R(\varepsilon)$,
$R>4$,
such  that
$$
u\leq \varepsilon
\quad
{\rm on}
\quad
\partial B(0,R).
$$
Set
$B=B(0,R)$,
$\widehat{B}\subset\subset B$.
Let
$v$
solve  the problem
$$
\left\{
\begin{array}{rcl}
\Delta v - v^q =0
&
{\rm in}
&
B\setminus K
\\
v(x)\to +\infty
&
{\rm when}
&
x \to K
\\
v=0
&
{\rm on}
&
\partial B.
\end{array}
\right.
$$
Then
$$
\Delta (v+\varepsilon)
-
(v+\varepsilon)^q
\leq 0 
\quad
{\rm in}
\quad
B \setminus K.
$$
Hence by asymptotic condition
(\ref{asympformaxsol})
and the comparison principle
\begin{equation}
\label{vapproxU}
u \leq v+\varepsilon
\quad
{\rm in}
\quad
B \setminus K.
\end{equation}
In what follows we first prove
(\ref{eta})
(\ref{upowerqtimeseta})
for
$v$
and then let
$\varepsilon$
vanish.

{\bf 3.}
Let
$\psi=1-\varphi$.
We claim that
\begin{equation}
\label{vpowerq}
\int_{B}
v^q \psi^m
\leq
C(m,n)\,
\capinRn (K)
\quad
{\rm for}
\quad
m\geq 2q'.
\end{equation}
In fact,
by Green's formula
\begin{eqnarray*}
\int_{B} v^q\psi^m
&
=
&
\int_{B}(\Delta v)\psi^m
\\
&
=
&
\int_{B}
v\Delta(\psi^m)
+
\int_{\partial B}
\left(
\psi^m
\frac{\partial v}{\partial \nu}
-
v
\frac{\partial \psi^m}{\partial \nu}
\right),
\\
\end{eqnarray*}
where 
$\nu$
is the outer normal on
$\partial B$.
Since
$\psi|_{\{ |x|\geq 2 \}}=1$
we conclude that
$$
\frac{\partial \psi^m}{\partial \nu}=0
\quad
{\rm on}
\quad
\partial B.
$$
By the comparison principle,
$v|_{B\setminus K}>0$.
Hence
$$
\int_{\partial B}
\psi^m
\frac{\partial v}{\partial \nu}\leq 0.
$$
Using the  H\"older inequality, we compute:
\begin{eqnarray}
\label{vpowerqintermid}
\int_{B}
v^{q}
\psi^m
&
\leq
&
\int_{B}
v\Delta (\psi^m)
\\
&
\leq
&
\int_{B}
v|\Delta (\psi^m)|
\nonumber
\\
&
\leq
&
m \int_{B}
\Big(
v\psi^{m-1}|\Delta\psi|
\Big)
+
m(m-1)\int_{B}
\Big(
v\psi^{m-2}| D \psi|^2
\Big)
\nonumber
\\
&
\leq
&
m
\left(
\int_{B}
v^{q}
\psi^m
\right)^{1/q}
\left(
\int_{\widehat{B}}
\psi^X
|\Delta \varphi|^{q'}
\right)^{1/q'}
\nonumber
\\
&
&
+
m(m-1)
\left(
\int_{B}
v^{q}
\psi^m
\right)^{1/q}
\left(
\int_{\widehat{B}}
\psi^Y
|D \varphi|^{2 q'}
\right)^{1/q'},
\end{eqnarray}
where 
$$
X=m-q',
\quad
Y=m-2q'.
$$
We can assume that
the left-hand side in
(\ref{vpowerq})
is positive.
From
(\ref{vpowerqintermid})
it then follows that
$$
\int_{B}
v^q \psi^m
\leq
m^{2q'}
\int_{\widehat{B}}
\left(
|\Delta \varphi|^{q'}
+
|D\varphi|^{2q'}
\right).
$$
Applying 
inequality
(\ref{interpolation}),
we obtain
$$
\int_{B}
v^q \psi^m
\leq
C(m,n)
\int_{\widehat{B}} |D^2 \varphi|^{q'},
$$
and 
(\ref{vpowerq})
follows from
(\ref{phi}).

{\bf 4.}
We claim that
\begin{equation}
\label{veta}
\int_{B}
v( |\Delta\eta| + |D\eta|)
\leq
C(m,n)
\capinRn (K).
\end{equation}
In fact, 
we have by the same calculations as in
(\ref{vpowerqintermid}):
\begin{eqnarray}
\int_{B} v|\Delta \eta|
&
\leq 
&
m
\left(
\int_{B}
v^{q}
\psi^{(m-1)q}
\right)^{1/q}
\left(
\int_{\widehat{B}}
|\Delta \varphi|^{q'}
\right)^{1/q'}
\nonumber
\\
& 
& 
+
m(m-1)
\left(
\int_{B}
v^{q}
\psi^{(m-2)q}
\right)^{1/q}
\nonumber
\\
&
&
\times
\left(
\int_{\widehat{B}}
|D \varphi|^{2q'}
\right)^{1/q'},
\label{1term}
\\
\int_{B}
v| D \eta|
&
\leq
&
m
\left(
\int_{B}
v^{q}
\psi^{(m-1)q}
\right)^{1/q}
\left(
\int_{\widehat{B}}
| D \varphi|^{q'}
\right)^{1/q'}.
\label{2term}
\end{eqnarray}
For
$m\geq 2q'$
we have
$$
(m-2)q\geq 2q',
\quad
(m-1)q\geq 2q'+q.
$$
Thus we can 
use
(\ref{vpowerq})
to
estimate the integrals containing
$v^q$
in
(\ref{1term})
and
(\ref{2term}).
Applying interpolation inequality
(\ref{interpolation})
to the last term in
(\ref{1term}),
we conclude on the basis of
(\ref{phi})
that
\begin{eqnarray*}
\int_{B}
v|\Delta\eta|
&
\leq
&
C(m,n)
\capinRn (K)^{1/q}
\left(
\int_{\widehat{B}}
|D^2\varphi|^{q'}
\right)^{1/q'}
\\
&
\leq
&
C(m,n)
\capinRn(K).
\end{eqnarray*}
Similarly,
applying the Poincar\'e 
inequality
to the last integral in
(\ref{2term})
gives
$$
\int_{B}
v|D\eta|
\leq
C(m,n)
\capinRn (K).
$$
We conclude that
(\ref{veta})
indeed holds.

{\bf 5.}
From
(\ref{vapproxU})
and
(\ref{veta})
we obtain
\begin{eqnarray*}
\int_{{\bf R}^n}
u (|D\eta| +|\Delta\eta|)
&
=
&
\int_{\widehat{B}} u (|D\eta| +|\Delta\eta|)
\\
&
\leq
&
C(m,n)
\left(
\capinRn (K)
+
\varepsilon 
\int_{\widehat{B}} (|D\eta| +|\Delta\eta|)
\right)
.
\end{eqnarray*}
To establish
(\ref{eta}) 
we
let
$\varepsilon\to 0$
both
in
(\ref{vapproxU})
and in the last inequality.
A similar limit argument applied to
(\ref{vpowerq})
gives us
(\ref{upowerqtimeseta}).
\end{proof}

Finally, 
we record a useful elementary inequality,
(see for example
\cite{Adams_Hedberg_book}
or
\cite{Maz'ya_book}).
Let
$J\in {\bf Z}$,
and let the function
$\zeta\colon (0, r_J)\to {\bf R}^1$
be either nondecreasing or
nonincreasing.
Then for any
$\kappa\in {\bf R}$
\begin{equation}
\label{sertointegr}
\sum_{j=J+1}^\infty
\zeta(r_j) 
r_j^{\kappa}
\lesssim
\int_0^{r_J}
\zeta (r) r^\kappa
\,
\frac{dr}{r}
\lesssim
\sum_{j=J}^\infty
\zeta (r_j) 
r_j^{\kappa}
.
\end{equation}

\section{Capacitary estimates}
\label{sec3}

\setcounter{equation}{0}

In this section we prove 
first estimates on
$u$
near
$\partial\Omega$.
We will work
in
$\mathbf{R}^n$
instead of
$\mathbf{S}^n$. 
According to  
sections~\ref{sphere&Rn},
\ref{capsection}
transition to the sphere is immmediate.

Theorems from this section 
will play the following role in 
the proof of the main result.
Let 
$u$
solve
\begin{equation}
\label{maineqagain}
u>0, 
\quad
\Delta u - u^q =0
\end{equation}
outside a compact set
$K\subset\mathbf{R}^n$.
When estimating the length of a 
curve
$\gamma$
in the metric
$u^{4/(n-2)}\geucl$
we will distinguish two regions.
In the first region
$\gamma$
is far enough from
$K$.
Then pointwise estimate 
(\ref{mainest})
from Theorem~\ref{captheorem}
will be applied.
In the second region
$\gamma $
is arbitrarily close
to
$K$.
Then we will use
integral estimate
(\ref{intu2/(n-2)est})
from
Theorem~\ref{intu2/(n-2)}.

\begin{theorem}
\label{captheorem}
Let
$K\subset B(0,r)$
be a compact set in
${\bf R}^n$,
$0<r<1$,
$n\geq 3$,
Let 
$u$
be the maximal solution 
of
(\ref{maineqagain})
in
$K^c$.
Then 
\begin{equation}
\label{mainest}
u(x)\asymp \frac{\capinRn(K)}{|x|^{n-2}},
\quad
|x|\geq 2r.
\end{equation}
\end{theorem}

\begin{proof}[of the upper estimate in (\ref{mainest})]
{\bf 1.}
According to the scalings
(\ref{scalingeq}),
(\ref{capscaling})
we need to prove 
that
for a compact set
$K$,
$K\subset B(0,1)$,
the following estimate holds:
\begin{equation}
\label{upper1}
u(x)
\lesssim
\capinRn (K)
\quad
{
\rm 
for
\quad
all
\quad
}
x
\quad
{\rm such \quad that}
\quad
2\leq |x| \leq 3.
\end{equation} 
Fix any such
$x$.
Let
$\eta$
be the function for our set
$K$
from
Lemma~\ref{mainlemma}
with some fixed
$m$.

{\bf 2.}
Utilising decay
(\ref{KellOsserm})
we can choose
$R>0$
so big that 
we have
\begin{eqnarray}
\label{upper2}
u(x)
&
=
&
(u\eta)(x)
\nonumber
\\
&
\lesssim
&
\int_{B(0,R)}
G_R(x,y)
\Delta 
(u\eta)(y)
\,dy.
\end{eqnarray}
Denote further
$B=B(0,R)$,
$G=G_R$.
Equation  
(\ref{maineqagain})
gives us 
\begin{eqnarray*}
\Delta(u \eta)
&
=
&
(\Delta u)\eta
+2\, Du \, D \eta
+u (\Delta\eta)
\\
&
\geq
&
2\, Du \, D\eta
+u\Delta\eta.
\end{eqnarray*}
Substitute this into
(\ref{upper2})
and integrate by parts to
deduce that
\begin{eqnarray*}
u(x)
&
\lesssim
&
-2
\int_B
D_y
G(x,y) 
D\eta(y)
u(y)
\,dy
\\
&
&
-
\int_B
G(x,y)
u(y)
\Delta\eta(y)
\,
dy.
\end{eqnarray*}
Next, the choice of 
$x$
and elementary bounds for
$G$
give that
\begin{equation*}
u(x)
\lesssim
\int_{\mathbf{R}^n}
u(|D \eta| +|\Delta\eta|).
\end{equation*}
Now estimate
(\ref{eta})
from
Lemma~\ref{mainlemma}
leads us to
(\ref{upper1}).
\end{proof}

\begin{proof}[of the lower estimate in (\ref{mainest})]
{\bf 1.}
According to the scalings
(\ref{scalingeq}),
(\ref{capscaling})
we need to prove 
that
for a compact set
$K$,
$K\subset B(0,1)$,
the following estimate holds:
\begin{equation}
\label{uinS}
u(x)
\gtrsim
\capinRn (K)
\quad
{
\rm 
for
\quad
all
\quad
}
x
\quad
{\rm such \quad that}
\quad
2\leq |x| \leq 3.
\end{equation} 
Taking  a suitable  approximation  
we can assume that
$K$
in
(\ref{uinS})
is the closure of a finite number of domains 
with smooth boundaries.

Now we recall the fundamental result in 
potential theory,
\cite{Adams_Hedberg_book} Ch. 2.
The Bessel kernel
$
{\mathcal J}_2\in 
C^\infty_\loc({\bf R}^n\setminus\{0\})
$
is defined 
via the formula
\begin{equation*}
(1-\Delta)^{-1} f
=
{\mathcal J}_2 * f
\quad
{\rm for}
\quad 
{\rm all}
\quad
f\in {\mathcal S}.
\end{equation*}
It satisfies the estimates
(see, for instance,
\cite{Adams_Hedberg_book} Chapter 1):
\begin{eqnarray}
\label{Besskern}
{\mathcal J}_2(x)
&
\asymp
&
{|x|^{-n+2}}
\quad
{\rm for}
\quad
x\in
B(0,1),
\\
{\mathcal J}_2(x)
&
\asymp
&
{e^{-|x|}} {|x|^{(-n+1)/2}}
\quad
{\rm for}
\quad
x\in
B(0,1)^c.
\nonumber
\end{eqnarray}
The theorem from nonlinear potential 
theory states that
there exists a Radon measure
$\mu^K$,
$\mu^K\geq 0$,
such that
\begin{equation*}
\supp(\mu^K)\subset K,
\end{equation*}
and
\begin{equation*}
\capinRn(K)
\asymp
\mu^K(K)
\asymp
\int_{\mathbf{R}^n}
\left(
\mathcal{J}_2*\mu^K
\right)^q.
\end{equation*}
Hence, after the regularisation
of
$\mu^K$
and a possible additional smooth approximation of
$K$
we obtain a function
$g\in C^\infty_0(\mathbf{R}^n)$,
$g\geq 0$,
such that
\begin{equation}
\label{gsupp}
\supp(g)\subset K,
\end{equation}
and
\begin{equation}
\label{gcap}
\capinRn(K)
\asymp
\int_{\mathbf{R}^n} g
\asymp
\int_{\mathbf{R}^n}
\left(
\mathcal{J}_2*g
\right)^q.
\end{equation}

{\bf 2.}
Set
$R=10$
and
$B=B(0,R)$.
For a fixed
$\varepsilon>0$
consider the Dirichlet problem
\begin{equation*}
\left\{
\begin{array}{rclll}
\Delta v 
&
=
& 
v^q -\varepsilon g
&
\
{\rm in}
\
&
B
\\
v
&
=
&
0
&
\
{\rm on}
\
&
\partial B.
\end{array}
\right.
\end{equation*}
As a simple consequence of the comparison principle
\cite{Marcus_Veron_Poincare_journal},
it has  the unique smooth solution
$v=v_\varepsilon$,
$v>0$
in
$B$.
Our goal will be to show that there exists
$\varepsilon>0$,
$\varepsilon =\varepsilon(n)$,
such that
\begin{equation}
\label{epsilonest}
v(x)\gtrsim \capinRn(K)
\quad
{\rm for\quad all}
\quad
x
\quad
{\rm such \quad that}
\quad
2\leq |x| \leq 3.
\end{equation}
To prove this 
we set
$G(x,y)=G_R(x,y)$,
and
note that 
by comparison principle
\begin{equation*}
v(x)
\leq
-\varepsilon
\,
\int_B
G(x,y)
g(y)
\,
dy
\quad
{\rm for\quad  all}
\quad
x\in B.
\end{equation*}
Consequently
\begin{eqnarray}
v(x)
&
=
&
-\varepsilon
\,
\int_B
G (x,y)
g(y)
\,
dy
+
\int_B
G (x,y)
v(y)^q
\,
dy
\nonumber
\\
&
\geq
&
\varepsilon
\,
\int_B
\big|
G (x,y)
\big|
\,
g(y)
\,
dy
\nonumber
\\
&
&
-
\varepsilon^q
\,
\int_B
\big|
G(x,y)
\big|
\,
\left(
\int_B
\big|
G(y,z)
\big|
\,
g(z)
\,
dz
\right)^q
\,
dy
\nonumber
\\
&
=
&
\varepsilon I(x)
-\varepsilon^q II(x)
\quad
{\rm for
\quad
all}
\quad
x\in B.
\label{epsilon,I,II}
\end{eqnarray}
Hence, 
to obtain
(\ref{epsilonest})
we need to estimate
$I$
from below
and
$II$
from above.

{\bf 3.}
Define
\begin{equation*}
S
=
\left\{
x\in\mathbf{R}^n:
\
2 \leq |x| \leq 3
\right\}.
\end{equation*}
The sets
$S$,
$\supp(g)$,
and
$\partial B$
are located at a distance at least
$1$
from each other.
Consequently
applying
(\ref{gsupp}),
(\ref{gcap}),
and invoking the elementary
properties
of 
$G$,
we derive 
\begin{eqnarray}
\label{I}
I(x)
&
=
&
\int_{\supp(g)}
\big|
G(x,y)
\big|
\,
g(y)
\, dy
\nonumber
\\
&
\gtrsim
&
\int_{\mathbf{R}^n}
g
\nonumber
\\
&
\gtrsim
&
\capinRn(K)
\quad
{\rm for \quad all}
\quad
x\in S.
\end{eqnarray}

{\bf 4.}
We claim that
\begin{equation}
\label{II}
II(x) \lesssim \capinRn (K)
\quad
{\rm for\quad all}
\quad
x\in S
.
\end{equation}
Indeed,  
fix
$x_0\in S$.
Introduce  the shell
\begin{equation*}
\widetilde{S}=
\left\{
x\in \mathbf{R}^n
:
\
2-1/100 \leq x \leq 3+ 1/100
\right\},
\end{equation*}
and utilise estimate
(\ref{Besskern})
for
$\mathcal{J}_2$
to write
\begin{eqnarray}
\label{II,X,Y}
II(x_0)
&
\lesssim
&
\int_{\widetilde{S}}
\frac{1}{|x_0 - y|^{n-2}}
\,
\left(
-
\int_B
G(y,z)g(z)
\, 
dz
\right)^q
\,
dy
\nonumber
\\
&
&
+
\int_{B\setminus\widetilde{S}}
\frac{1}{|x_0 - y|^{n-2}}
\,
\big(
\mathcal{J}_2 * g
\big)^q
(y)
\,
dy
\nonumber
\\
&
=
&
X+Y.
\end{eqnarray}
We estimate
$X$
and
$Y$
separately.

To estimate 
$X$
define
the function
$H\colon B\to \mathbf{R}^1$
by writing
\begin{equation*}
H(y)=-\int_B
G(y,z)g(z)
\, 
dz,
\quad
y\in B.
\end{equation*}
Notice that according to
(\ref{gsupp}), 
$H$
is positive 
in
$B$
and harmonic in
$B\setminus K$.
Consequently
$H^q$
is subharmonic
in
$B\setminus K$.
Hence by the mean value property
\begin{eqnarray*}
X
&
\lesssim
&
\left(
\max_{\widetilde{S}} H
\right)^q
\int_{\widetilde{S}}
\frac{dy}{|x_0-y|^{n-2}}
\\
&
\lesssim
&
\max_{\widetilde{S}}
H^q
\\
&
\lesssim
&
\int_B
H(y)^q\,dy.
\end{eqnarray*}
Now 
(\ref{Besskern})
and 
(\ref{gcap})
allow  us  to conclude that
\begin{eqnarray}
\label{X}
X
&
\lesssim
&
\int_{\mathbf{R}^n}
\big(
\mathcal{J}_2 * g
\big)^q(y)
\,
dy
\nonumber
\\
&
\lesssim
&
\capinRn(K).
\end{eqnarray}

To estimate
$Y$
notice that
\begin{equation*}
|x_0 -y |\geq 1/100
\quad
{\rm for\quad all}
\quad
y\in B\setminus\widetilde{S}.
\end{equation*}
Therefore
utilising
(\ref{gcap})
we derive
\begin{eqnarray}
\label{Y}
Y
&
\lesssim
&
\int_{B\setminus\widetilde{S}}
\big(
\mathcal{J}_2 * g
\big)^q(y)
\,
dy
\nonumber
\\
&
\lesssim
&
\int_{\mathbf{R}^n}
\big(
\mathcal{J}_2 * g
\big)^q
\nonumber
\\
&
\lesssim
&
\capinRn(K).
\end{eqnarray}
Substituting
(\ref{X})
and
(\ref{Y})
into
(\ref{II,X,Y})
we deduce
(\ref{II}).

{\bf 5.}
Now we conclude the proof of the theorem.
First we  establish
(\ref{epsilonest}).
Substitute
(\ref{I})
and
(\ref{II})
into
(\ref{epsilon,I,II}):
\begin{equation*}
v(x)
\geq
\big(
\varepsilon C_1(n)- \varepsilon^q C_2(n)
\big)
\capinRn(K)
\quad
{\rm for \quad all}
\quad
x\in S.
\end{equation*}
Choosing
the suitable
$\varepsilon>0$
derive
(\ref{epsilonest}).

Finally,
the regularity of
$K$
implies
that
our maximal solution
$u$
blows up near 
$K$
as in
(\ref{asympformaxsol}).
Therefore
\begin{equation*}
u\geq v
\quad
{\rm on}
\quad
\partial (B\setminus K).
\end{equation*}
Owing to 
(\ref{gsupp})
and the comparison principle,
\begin{equation*}
u\geq v
\quad
{\rm in}
\quad
B\setminus K.
\end{equation*}
This inequality and
(\ref{epsilonest})
complete the proof of
(\ref{uinS}).
\end{proof}

Next we establish an integral estimate for any 
solution of
(\ref{maineqagain}).
This is
Theorem~\ref{intu2/(n-2)}
below.  It has a particularly simple proof when
$n\geq 4$
and hence
\begin{equation*}
\frac{2}{n-2}\leq 1.
\end{equation*}
In this case it follows 
more or less
directly from 
the representation formula for
solution of the linear Poisson equation.
However, such  approach does not work for
$n=3$
because
the singularity of the Green function is 
too strong then.
Proof of 
Theorem~\ref{intu2/(n-2)}
given below does not use representation 
formula. Instead we rely on techniques
common in 
quasilinear elliptic regularity theory. Such arguments  
were first used by
Moser
\cite{Moser_1},
\cite{Moser_2}
for linear equations, and by
Trudinger
\cite{Trudinger_weak_Harnack}
for nonlinear equations.

\begin{theorem}
\label{intu2/(n-2)}
Let
$K\subset B(0,1)$
be a compact set,
$\varphi$
be the function from
Lemma~\ref{mainlemma},
and let
\begin{equation*}
m=\frac{n+2}{2} +100n.
\end{equation*}
Then for any 
$u$
solving 
(\ref{maineqagain})
in
$K^c$
the estimate
\begin{equation}
\label{intu2/(n-2)est}
\int_{B(0,10)}
u^{2/(n-2)}
\,
(1-\varphi)^m
\lesssim
\capinRn(K)^{2/(n-2)}
\end{equation}
holds.
\end{theorem}

\begin{proof}
{\bf 1.}
We claim that for any number
$\varepsilon$,
$0<\varepsilon<1$,
the inequality
\begin{eqnarray}
\label{iterstep}
\left(
\int_{\mathbf{R}^n}
u^{(1-\varepsilon)n/(n-2)}
\,
|\zeta|^{2n/(n-2)}
\right)^{(n-2)/n}
&
\leq
&
C(\varepsilon)
\Big(
\int_{\mathbf{R}^n}
u^{1-\varepsilon}
|D\zeta|^2
\nonumber
\\
&
&
+
\int_{\mathbf{R}^n}
u^{q-\varepsilon}
\zeta^2
\Big)
\end{eqnarray}
holds for all functions
$\zeta\in C^\infty_0 ( \mathbf{R}^n )$,
such that
$\zeta=0$
in an open neighbourhood of
$K$.
In fact,
multiplying
the equation
\begin{equation*}
\Delta u - u^q = 0
\quad
{\rm in}
\quad
\mathbf{R}^n\setminus K 
\end{equation*}
by
$u^{-\varepsilon}\zeta^2$,
integrating by parts, and 
invoking the formula
\begin{equation*}
D(u^{-\varepsilon}\zeta^2)
=-\varepsilon u^{-\varepsilon-1}
\zeta^2
Du
+
2\zeta u^{-\varepsilon} D\zeta,
\end{equation*}
we deduce that
\begin{eqnarray*}
\varepsilon
\int_{\mathbf{R}^n}
|Du|^2
u^{-\varepsilon-1}
\zeta^2
&
\leq 
&
2
\int_{\mathbf{R}^n}
|Du|
\,
|D\zeta|
\,
u^{-\varepsilon}
\,
|\zeta|
\\
&
&
+
\int_{\mathbf{R}^n}
u^{q-\varepsilon}
\zeta^2.
\end{eqnarray*}
Since
\begin{equation*}
|Du|
\,
|D\zeta|
\,
u^{-\varepsilon}
\,
|\zeta|
\leq
\delta
|Du|^2
u^{-\varepsilon-1}
\zeta^2
+
\frac{1}{4\delta}
|D\zeta|^2
u^{-\varepsilon+1}
\end{equation*}
for each
$\delta>0$,
we derive that
\begin{equation*}
\int_{\mathbf{R}^n}
|Du|^2
\,
u^{-\varepsilon -1}
\zeta^2
\leq
C(\varepsilon)
\left(
\int_{\mathbf{R}^n}
u^{1-\varepsilon}
|D\zeta|^2
+
\int_{\mathbf{R}^n}
u^{q-\varepsilon}
\zeta^2
\right).
\end{equation*}
After some calculations we find
\begin{equation*}
\int_{\mathbf{R}^n}
\left|
D
\left(
u^{(1-\varepsilon)/2}
\,
\zeta
\right)
\right|^2
\leq
C(\varepsilon)
\left(
\int_{\mathbf{R}^n}
u^{1-\varepsilon}
|D\zeta|^2
+
\int_{\mathbf{R}^n}
u^{q-\varepsilon}
\zeta^2
\right).
\end{equation*}
Now
(\ref{iterstep})
follows from the Sobolev
inequality
applied to the
left hand side.

{\bf 2.}
Set
$B=B(0,10)$
and
$\widehat{B}=B(0,20)$.
In
(\ref{iterstep})
choose
$\varepsilon \in (0,1)$
such that
\begin{equation*}
(1-\varepsilon)
\frac{n}{n-2}
=\frac{2}{n-2}.
\end{equation*}
Then select
a smooth cutoff function
$\theta\in C^\infty_0(\widehat{B})$,
such that
$\theta = 1$
on
$B$.
Take the function
$\eta=(1-\varphi)^m$
from
Lemma~\ref{mainlemma}.
Now 
set
$\zeta = \eta\theta$,
in estimate
(\ref{iterstep})
to discover that
\begin{eqnarray}
\label{almostiterstep}
\left(
\int_{B}
u^{2/(n-2)}
\,
(1-\varphi)^m
\right)^{(n-2)/2}
\lesssim
&
\Big(
&
\int_{\mathbf{R}^n}
u^{1-\varepsilon}
|D\eta|^2
\nonumber
\\
&
&
+
\int_{\mathbf{R}^n}
u^{1-\varepsilon}
|D\theta|^2
\\
&
&
+
\int_{\mathbf{R}^n}
u^{q-\varepsilon}
(\eta\theta)^2
\
\Big)
^{1/(1-\varepsilon)}.
\nonumber
\end{eqnarray}
We estimate three integrals 
in the right hand side of
(\ref{almostiterstep})
as follows.
Applying Holder inequality  and estimates
(\ref{phi}), 
(\ref{upowerqtimeseta})
from
Lemma~\ref{mainlemma}
we deduce that
\begin{eqnarray*}
\int_{\mathbf{R}^n}
u^{1-\varepsilon}
|D\eta|^2
&
\lesssim
&
\left(
\int_B u^{q(1-\varepsilon)}
(1- \varphi)^{(n+2)/2}
\right)^{1/q}
\left(
\int_B
|D\varphi|^{2q'}
\right)^{1/q'}
\\
&
\lesssim
&
\capinRn (K)^{(q-\varepsilon)/q}.
\end{eqnarray*}
Estimate 
(\ref{mainest})
from
Theorem~\ref{captheorem}
implies
\begin{eqnarray*}
\int_{\mathbf{R}^n}
u^{1-\varepsilon}
|D\theta|^2
&
\lesssim
&
\| 
u
\|^{1-\varepsilon}_{L^\infty(\widehat{B}\setminus B)}
\\
&
\lesssim
&
\capinRn (K)^{1-\varepsilon}.
\end{eqnarray*}
Finally, 
owing to Holder inequality and
(\ref{upowerqtimeseta})
we have
\begin{equation*}
\int_{\mathbf{R}^n}
u^{q-\varepsilon}
\,
(\eta\theta)^2
\lesssim
\capinRn (K)^{(q-\varepsilon)/q}.
\end{equation*}
Substituting  these estimates in
(\ref{almostiterstep})
we arrive at
\begin{eqnarray*}
\left(
\int_{B}
u^{2/(n-2)}
\,
(1-\varphi)^m
\right)^{(n-2)/2}
&
\lesssim
&
\capinRn (K)
+
\capinRn (K)^{(q-\varepsilon)/(q-\varepsilon q)}
\\
&
\lesssim
&
\capinRn (K).
\end{eqnarray*}
This is assertion
(\ref{intu2/(n-2)est}).
\end{proof}

\section{Proof of the Wiener test for conformal metrics}

%
%
%
%

\setcounter{equation}{0}

\subsection{Sufficiency}

We prove the implication
$(ii)\Rightarrow(i)$
in Theorem~\ref{maintheorem}.

{\bf 1.}
In the proof we will work on
$\mathbf{S}^n$.
A curve
$\gamma\colon[0,+\infty)\to \Omega$
is said to converge to infinity if for every
compact set
$M\subset\Omega$,
there is a time
$T$,
$0<T<+\infty$,
such that
$\gamma(t)\not\in M$
for all
$t>T$.
By a version of the Hopf-Rinow theorem,
a metric is complete in
$\Omega$
if and only if
every smooth curve
converging to infinity
has the infinite length.

Let 
$u_{\Omega}$ 
be the maximal solution
of the conformal scalar curvature equation
(\ref{eq_on_sphere})
in 
$\Omega$.
We set
\begin{equation*}
g=u_{\Omega}^{4/(n-2)}\gstand  
\end{equation*}
Fix any smooth curve
$\gamma\colon[0,+\infty)\to \Omega$
converging to infinity.
To prove statement 
(i) in   
Theorem~\ref{maintheorem}
we need to show that
\begin{equation}
\label{lengthinf1}
L_g(\gamma)=
\int_0^\infty
u_{\Omega}(\gamma)^{2/(n-2)}
\,
\gstand(\dot{\gamma}, \dot{\gamma})^{1/2}
\,
dt
=
+\infty.
\end{equation}
In the rest of the proof we establish
(\ref{lengthinf1}).

{\bf 2.}
Compactness of
$\mathbf{S}^n$ 
and convergence of
$\gamma$
to infinity imply
the existence of a point
$p\in K$
such that
\begin{equation}
\label{distinf}
d_{\gstand}(\gamma (T_k), p)\to 0
\quad
{\rm for\quad a \quad sequence}
\quad
\{T_k\},
\quad
T_k\to+\infty.
\end{equation}
For
$j=1$,
$2$,
$\ldots$
we define
$\Gamma_j$
to be that part of
$\gamma$
whose image is contained in
the shell
$S_j$,
\begin{equation*}
S_j=
\left\{
x\in {\bf S}^n:
r_j<d_{\gstand}(x,p)<r_{j-1}
\right\}.
\end{equation*}
The smoothness of
$\gamma$
implies that
for any 
$j\geq j_0$
the set
$\Gamma_j$
is at most a countable union
of open smooth curves. 
Utilising condition
(\ref{distinf})
we deduce   
that
$\Gamma_j\ne\emptyset$,
and moreover
\begin{equation*}
L_{\gstand}(\Gamma_j)\geq \frac{r_j}{100}
\quad
{\rm
for 
\quad
all 
}
\quad
j \geq j_0.
\end{equation*}
After a rotation we 
can assume
that
\begin{equation*}
K\cap \overline{B}(p, r_j)\subset \mathbf{U}
\quad
{\rm for \quad all \quad}
 j\geq j_0-10,
\end{equation*}
where
$\mathbf{U}$
is cup
(\ref{spherical_cup})
around the south pole from the definition of
$\capacity$.
We claim that 
for 
all 
$j \geq j_0$
the inequality
\begin{equation}
\label{lengthj}
L_{g}(\Gamma_j)
\gtrsim
r_j
\left(
\frac
{
\capacity
\left( 
K  \cap \overline{B}(p, r_{j+2}) 
\right)
}
{
r_j^{n-2}
}
\right)
^{2/(n-2)}
\end{equation}
holds.
In fact,
define
the open set
$\Omega_j$,
$\Omega_j\supset\Omega$,
by writing
\begin{equation*}
\Omega_j=
{\bf S}^n\setminus
( K \cap \overline{B}(p,r_{j+2}) )
.
\end{equation*}
Let
$u_j$
be the maximal solution to our equation
(\ref{eq_on_sphere})
in 
$\Omega_j$.
Pull
estimate
(\ref{mainest})
from Theorem~\ref{captheorem}
back
to the sphere
via the
stereographic projection,
keeping in mind  that the conformal factor
in
(\ref{conf_change_of_conf_factor})
satisfies
\begin{equation*}
\Upsilon(x)\asymp 1
\quad
{\rm for \quad all \quad }
x,
\quad
|x|\leq 10.
\end{equation*}
We discover  that
\begin{eqnarray*}
u_{\Omega}(x)
&
\geq
&
u_j(x)
\\
&
\gtrsim
&
\frac
{\capacity
\left(K \cap
\overline{B}(p, r_{j+2}) 
\right)
}
{
r_j^{n-2}
}
\quad
{\rm for \quad all}
\quad
x\in S_j \cap \Omega.
\end{eqnarray*}
Let 
$I_j$,
$I_j\subset (0,+\infty)$,
be the open set such that
\begin{equation*}
\Gamma_j\colon I_j\to \Omega\cap S_j.
\end{equation*}
Then we derive that
\begin{eqnarray*}
L_{g}(\Gamma_j)
&=&
\int_{I_j}
u_{\Omega}(\gamma)^{2/(n-2)}
\,
\gstand(\dot{{\gamma}},\dot{{\gamma}} )^{1/2}
\,
dt  
\\
&\geq& 
\left(
\inf_{S_j} u_{\Omega} 
\right)^{2/(n-2)}
L_{\gstand}(\Gamma_j)
\\
&\gtrsim&
\left(
\frac
{\capacity
\left( 
K \cap \overline{B}(p, r_{j+2}) 
\right)
}
{
r_j^{n-2}
}
\right)^{2/(n-2)}
\, 
r_j,
\end{eqnarray*}
thereby obtaining
(\ref{lengthj}).

{\bf 3.}
We claim that
(\ref{lengthinf1})
holds.
Indeed, 
the sets
$S_j$
are disjoint, 
and thus
\begin{equation*}
L_{g}({\gamma})
\geq
\sum_{j\geq 1}
L_{g}(\Gamma_j).
\end{equation*}
To each term with sufficiently large number  
in this sum we 
apply estimate 
(\ref{lengthj})
and recall
(\ref{capball}) 
to  derive that
\begin{eqnarray*}
L_{g}(\gamma) 
&\gtrsim&
\sum_{j\geq j_0}
r_{j+2}
\left(
\frac
{
\capacity
\left( 
K \cap \overline{B}(p, r_{j+2}) 
\right)
}
{
r_{j+2}^{n-2}
}
\right)^{2/(n-2)}
\\
&\gtrsim&
\sum_{j\geq j_0+100}
\left(
\frac
{\capacity(B(p, r_j) \cap K )
}
{
\capacity(B(p,r_j))
}
\right)^{2/(n-2)}.
\end{eqnarray*}
Finally utilise
(\ref{sertointegr})
and  
(\ref{wienertest})
to establish
(\ref{lengthinf1}).
This completes the proof of 
implication
$(ii)\Rightarrow (i)$
in Theorem~\ref{maintheorem}.

\subsection{Necessity}

Now we prove the implication
$(i)\Rightarrow (ii)$
in
Theorem~\ref{maintheorem}.

{\bf 1.}
Seeking 
a contradiction assume that
$(ii)$
does not hold.
Hence
\begin{equation}
\label{negatwiener}
\int_0^{1/2}
\left(
\frac{\capacity(B(P, r)\cap K)}{\capacity(B(P,r))}
\right)^{2/(n-2)}
\,
\frac{dr}{r} < +\infty
\end{equation}
for
some
$P\in\partial\Omega$.
The desired contradiction will follow if 
the maximal solution
of 
(\ref{eq_on_sphere})
does not give the metric complete in
$\Omega$. 
Let
$U$
be this maximal solution, 
and let
\begin{equation*}
g=U^{2/(n-2)}\gstand .
\end{equation*}
According to 
the Hopf-Rinow theorem,
to prove the non-completeness of
$g$
we must show that
there exists a smooth curve
$c$,
\begin{eqnarray}
c\colon [0,1)\to \Omega,
&
\
&
{\rm 
such
\quad 
that
\quad
}
d_{\gstand}(c(t), P)\to 0 
{
\quad
\rm 
as
\quad
}
t\to 1,
\nonumber
\\
&
&
{\rm 
and
\quad
}
L_g(c)<+\infty.
\label{finlengthatP}
\end{eqnarray}

{\bf 2.}
First we reformulate
claim
(\ref{finlengthatP}).
Fix a parameter
$\rho>0$,
which we will later  choose small.
Set
\begin{equation*}
\widetilde{K}=K\cap\overline{B}(P,\rho).
\end{equation*}
Let
$U_1$
be the maximal solution of 
(\ref{eq_on_sphere})
in
$
\widetilde{\Omega}
$,
\begin{equation*}
\widetilde{\Omega}
=
\mathbf{S}^n 
\setminus 
\widetilde{K},
\end{equation*}
and let
$U_2$
be the maximal solution of
(\ref{eq_on_sphere})
in
$\Omega \cup B(P,\rho)$.
From
(\ref{sum})
we deduce that
\begin{equation*}
U^\metricexp
\lesssim
U_1^\metricexp
+
U_2^\metricexp
\quad
{\rm in}
\quad
\Omega.
\end{equation*}
At the same time 
(\ref{KellOsserm}) 
implies
\begin{equation*}
U_2(x)
\leq 
C(\rho)
\quad
{\rm 
for
\quad 
all}
\quad
x
\in
B(P, \rho / 2)
.
\end{equation*}
Therefore
(\ref{finlengthatP})
is equivalent to the same 
statement with
$\Omega$ 
replaced by
$\widetilde{\Omega}$,
and 
$g$
replaced
by
\begin{equation*}
\widetilde{g}
=U^\metricexp_1
\gstand.
\end{equation*}
To prove this statement it will be convinient 
to  transform  the problem to
${\bf R}^n$.

Applying a  suitable rotation and stereographic projection we can achieve that
$P$
is mapped to
$0$.
We denote the image of
$K$
under such
map by the same letter
$K$.
In
$\mathbf{R}^n$
we set
\begin{equation*}
K_j=K\cap \overline{B}_j
\quad
{\rm and}
\quad
B=B(0,1).
\end{equation*}
By 
$u$
we 
denote  the conformal pullback  
(\ref{conf_change_of_conf_factor})
of
$U_1$,
\begin{equation*}
u(x)
=
\Upsilon(x)
U_1 (\sigma^{-1}x),
\quad
x
\in
\sigma(\widetilde{\Omega}).
\end{equation*}
As it is shown in 
section~\ref{sphere&Rn},
$u$
is the maximal solution of 
(\ref{positsoltomaineq})
in
$\mathbf{R}^n \setminus K_J$
for some
$J$,
$J=J(\rho)$.
From
(\ref{negatwiener})
and 
(\ref{capequivalence})
we deduce that
\begin{equation}
\label{thinKJ}
\int_0^1
\left(
\frac
{\capinRn(B(0,r)\cap K_J)}
{\capinRn(B(0,r))}
\right)^\metricexp
\,
\frac{dr}{r}
<+\infty.
\end{equation}
Finally, to establish
(\ref{finlengthatP})
we must prove that
$u^\metricexp\geucl$
is not complete, that is
\begin{eqnarray}
\label{contr:finlength}
\int_\gamma
u^\metricexp
\, ds
<+\infty
&
&
{\rm for \ a \ smooth \ curve \ } 
\gamma\colon [0,1)\to B\setminus K_J,
\nonumber
\\
&
&
{\rm such\  that \ }
\gamma(t)\to 0
\ 
{\rm as}
\ 
t\to 1.
\end{eqnarray}

{\bf 3.}
We intend to establish
(\ref{contr:finlength}).
The construction
of 
$\gamma$
in
(\ref{contr:finlength})
will be indirect. More precisely,  
let us first  reduce the proof of
(\ref{contr:finlength})
to an integral estimate for
our maximal solution
$u$.

We 
{\it assert} 
that it is possible to choose large enough
$J$
in
(\ref{thinKJ})
(equivalently, to choose small enough
$\rho>0$)
such that there exists a compact set
$\Sigma$,
\begin{equation*}
K_J\subset\Sigma\subset B,
\end{equation*}
with the following two properties:
\begin{equation}
\int_{B\setminus\Sigma}
u(x)^\metricexp
\,
\frac{1}{|x|^{n-1}}
\, dx
<+\infty,
\label{Sigma1}
\end{equation}
and
\begin{equation}
\Hmeas^{n-1}
\big(
\pi
(
\Sigma\setminus\{0\}
)
\big)
<
\Hmeas^{n-1}(\partial B),
\label{Sigma2}
\end{equation}
where
$\pi$
is the radial projection on
$\partial B$,
\begin{equation*}
\pi\colon
B\setminus\{0\}\to \partial B,
\quad
x\mapsto\frac{x}{|x|}.
\end{equation*}
This assertion is the core of the proof.
Before passing to its verification
we conclude the current step by showing
that
(\ref{Sigma1}),
(\ref{Sigma2})
immediately imply 
(\ref{contr:finlength})
and hence the theorem.

Indeed, for 
$\omega\in\partial B$
we define the interval
$\ell(\omega)$
by writing
\begin{equation*}
\ell(\omega)=
\left\{
x\in {\bf R}^n:
\quad
x=s\omega, 
\quad
0<s\leq 1
\right\}.
\end{equation*}
Set
\begin{equation*}
\Xi=
\partial B 
\setminus 
\pi(\Sigma \setminus\{0\}).
\end{equation*}
First notice that
\begin{equation*}
\pi^{-1}
( \Xi ) 
\subset 
B\setminus\Sigma.
\end{equation*}
Hence,
using the  polar coordinates 
$(r,\omega)$,
$r>0$,
$\omega\in\partial B$,
we deduce at once 
from
(\ref{Sigma1})
that
\begin{eqnarray*}
+\infty
&
>
&
\int_{\pi^{-1}( \Xi ) }
u(x)^\metricexp
\,
\frac{1}{|x|^{n-1}}
\, dx
\\
&
=
&
\int_{\Xi}
\int_0^1
u(x(r,\omega))^\metricexp
\,
\frac{1}{r^{n-1}}
\,
r^{n-1}
\, 
dr
\, 
d\Hmeas^{n-1}(\omega)
\\
&
=
&
\int_{\Xi}
\left(
\int_{\ell(\omega)}
u^\metricexp
\, 
ds
\right)
\, d\Hmeas^{n-1}(\omega).
\end{eqnarray*}
Next, apply
(\ref{Sigma2})
to discover that
\begin{equation*}
\Hmeas^{n-1} (\Xi) 
=
\Hmeas^{n-1}(\partial B)
-
\Hmeas^{n-1}
\big( 
\pi \left(
\Sigma\setminus\{0\}
\right) 
\big)
> 0.
\end{equation*}
Consequently
\begin{equation*}
\int_{\ell(\omega_0)}u^\metricexp
\, ds
<+\infty
\quad{\rm for\quad some\quad}
\omega_0 \in \Xi.
\end{equation*}
By our definitions
\begin{equation*}
\ell(\omega_0)\cap K_J =\emptyset,
\end{equation*}
and we conclude that 
(\ref{contr:finlength})
holds for the curve
$\gamma=\ell(\omega_0)$.

Thus, to establish  the theorem it is left to 
construct
$\Sigma$
satisfying
(\ref{Sigma1}) 
and
(\ref{Sigma2}).
The rest of the proof is devoted entirely to this
construction.

{\bf 4.}
For any 
$j\geq J$
let
\begin{equation*}
v=v_j
\end{equation*}
be the maximal solution for
$K_{j-2}$.
We now  apply
Theorem~\ref{intu2/(n-2)}
and scalings
(\ref{scalingeq}),
(\ref{capscaling})
to bound
$v$.
Applying
the scaling
\begin{equation*}
x\mapsto r_{j-3}\,x,
\quad
x\in \mathbf{R}^n,
\end{equation*}
to estimate
(\ref{intu2/(n-2)est})
we deduce 
that
there exists a function
\begin{equation*}
\varphi=\varphi_j,
\end{equation*}
such that:
\begin{eqnarray}
&&
{\varphi \in C^\infty_0(B(0,2))}, 
\quad
0\leq\varphi\leq 1
\quad
{\rm in}
\quad
B(0,2),
\nonumber
\\
&&
\varphi =1
\quad
{\rm in}
\quad
{\rm a\quad  neighbourhood\quad  of\quad  }
K_{j-2},
\nonumber
\\
&&
\int_{B(0,2)}
\left|
D^2\varphi
\right|
^
{(n+2)/4}
\lesssim
\capinRn (K_{j-2}),
\label{choiceofphi1}
\end{eqnarray}
and
\begin{eqnarray}
\frac{1}{r_j^{n-1}}
\,
\int_{B_{j-2}}
v^\metricexp
(1-\varphi)^m
&
\lesssim
&
\big(
\capinRn (K_{j-2}/r_{j-3})
\big)^\metricexp
\nonumber
\\
&
\lesssim
&
\left(
\frac{\capinRn (K_{j-2})}{\capinRn (B_{j-2})}
\right)^\metricexp
.
\label{choiceofphi2}
\end{eqnarray}

{\bf 5.}
Now we construct the compactum
$\Sigma$
for
(\ref{Sigma1}),
(\ref{Sigma2}).
Set
\begin{equation*}
S_j=
\left\{
x:
\
r_j\leq |x| \leq r_{j-1} 
\right\}.
\end{equation*}
First for a fixed
$j\geq J$,
define the
compact set
$E_j$
by writing
\begin{equation*}
E_j
=
\left\{
x\in S_j:
\
\varphi_j(x)\geq 
\frac
{99}
{100}
\right\},
\end{equation*}
where the function
$\varphi_j$
is taken from
(\ref{choiceofphi1}),
(\ref{choiceofphi2}).
Then define
\begin{equation*}
\Sigma
=
\left(
\bigcup_{j=J}^\infty E_j
\right)
\bigcup
\{0\}.
\end{equation*}
According to the construction,
the set
$\Sigma$
is compact  and
\begin{equation*}
K_J\subset\Sigma\subset B.
\end{equation*}
We claim that
(\ref{Sigma2})
holds.
Indeed, the definition of the capacity
and
(\ref{choiceofphi1})
imply that
\begin{eqnarray*}
\capinRn (E_j)
&
\leq
&
\int_{B(0,2)}
\left|
D^2
\left(
\frac
{100}
{99}
\,
\varphi_j
\right)
\right|^{(n+2)/4}
\\
&
\lesssim
&
\capinRn (K_{j-2}).
\end{eqnarray*}
The metric estimate
(\ref{capHausdcontent})
therefore ensures
\begin{equation*}
\Hmeas^{n-1}_{\infty}(E_j)^{(n-2)/2}
\lesssim
\capinRn (K_{j-2})^{n-1}.
\end{equation*}
The projection
$\pi$
restricted to
$S_j$
distorts the distances
at most
$1/r_j$
times.
Consequently
\begin{eqnarray*}
\Hmeas^{n-1}_{\infty}
(
\pi (\Sigma\setminus\{0\})
)
&
\leq
&
\sum_{j=J}^\infty
\Hmeas^{n-1}_{\infty}(\pi(E_j))
\\
&
\leq
&
\sum_{j=J}^\infty
\frac{1}{r_j^{n-1}}
\,
\Hmeas^{n-1}_{\infty}(E_j)
\\
&
\lesssim
&
\sum_{j= J}^\infty
\frac{1}{r_j^{n-1}}
\,
\capinRn
(K_{j-2})
^{(n-1)\metricexp}
\\
&
\lesssim
&
\left(
\sum_{j=J-2}^\infty
\left(
\frac
{
\capinRn(K_j)
}
{
\capinRn
(B_j)
}
\right)
^{\metricexp}
\right)^{n-1}
.
\end{eqnarray*}
According to 
(\ref{sertointegr})
and
(\ref{thinKJ})
we can make 
the last series
as small as we wish
by choosing
$J$
large enough.
Thus  for any
$\varepsilon>0$
we may fix
$J$
in
(\ref{thinKJ})
so that
\begin{equation*}
\Hmeas^{n-1}_{\infty}
\big(
\pi 
(
\Sigma\setminus \{0\}
)
\big)<\varepsilon.
\end{equation*}
For sets lying on an
$s$-dimensional smooth submanifold,
the Hausdorff
$s$-measure is
equivalent to the Lebesgue 
$s$-measure.
Consequently
\begin{equation*}
\Hmeas^{n-1}_{\infty}(F)
\gtrsim
\Hmeas^{n-1}(F)
\quad
{\rm for\quad any}
\quad
F\subset\partial B.
\end{equation*}
This gives
(\ref{Sigma2}).

{\bf 6.}
It is left to prove
(\ref{Sigma1}).
Splitting  the integral there  we find that 
\begin{equation}
\label{Sigma1proof1}
\int_{B\setminus\Sigma}
u(x)^\metricexp
\,
\frac{1}{|x|^{n-1}}
\, dx
\lesssim
\sum_{j=J}^\infty
\frac{1}{r_j^{n-1}}
\,
\int_{S_j\setminus\Sigma}
u^\metricexp
.
\end{equation}
Thus our task is to estimate
$u$
in
$S_j\setminus\Sigma$.
Fix 
$j\geq J$.
Let as before 
$v$
be the maximal solution for
$K_{j-2}$.
For
$l=1$,
$2$,
$\ldots$,
$j-2$
let
$w_l$
be the maximal solution for
$K\cap S_l$.
From 
(\ref{sum})
we deduce that
\begin{equation*}
u\leq v
+
\sum_{l=1}^{j-2} w_l
\quad
{\rm in}
\quad
S_j.
\end{equation*}
Next observe  that
\begin{equation*}
|x-y| \asymp r_l
\quad
{\rm for
\quad 
all}
\quad
x\in S_j,
\quad
y\in S_l,
\quad
l\leq j-2.
\end{equation*}
Hence applying the scaled estimate
(\ref{mainest}) 
from
Theorem~\ref{captheorem}
to
$w_l$,
we derive that
\begin{eqnarray}
\label{shellinfluence}
u(x)
&
\lesssim
&
v(x) +\sum_{l=1}^{j-2}
\frac{\capinRn (K\cap S_l)}
{r_l^{n-2}}
\nonumber
\\
&
\lesssim
&
v(x)
+
\sum_{l=1}^{j-2}
\frac{\capinRn (K_l)}{r_l^{n-2}}
\quad
{\rm for \quad all}
\quad
x\in S_j.
\end{eqnarray}
To estimate
$v$
in
$S_j$
notice that
(\ref{choiceofphi2})
and the definition of
$\Sigma$
ensure that
\begin{equation*}
\frac{1}{r_j^{n-1}}
\,
\int_{S_j\setminus\Sigma}
v^\metricexp
\lesssim
\left(
\frac
{
\capinRn(K_{j-2})
}
{
\capinRn(B_{j-2})
}
\right)^\metricexp
.
\end{equation*}
Utilising
(\ref{shellinfluence})
we thereupon conclude that
\begin{eqnarray*}
\frac{1}{r_j^{n-1}}
\,
\int_{S_j\setminus\Sigma}
u^\metricexp
&
\lesssim
&
\frac{1}{r_j^{n-1}}
\,
\int_{S_j\setminus\Sigma}
v^\metricexp
\\
&
&
+
\frac{1}{r_j^{n-1}}
\left(
\sum_{l=1}^{j-2}
\frac{\capinRn (K_l)}{r_l^{n-2}}
\right)^\metricexp
\,
|S_j|
\\
&
\lesssim
&
\left(
\frac
{
\capinRn(K_{j-2})
}
{
\capinRn(B_{j-2})
}
\right)^\metricexp
\\
&
&
+
r_j
\left(
\sum_{l=1}^{j-2}
\frac{\capinRn (K_l)}{r_l^{n-2}}
\right)^\metricexp
.
\end{eqnarray*}
Now continue
(\ref{Sigma1proof1})
to find
\begin{eqnarray}
\label{Sigma1proof2}
\int_{B\setminus \Sigma}
u(x)^\metricexp
\,
\frac{1}{|x|^{n-1}}
\, 
dx
&
\lesssim
&
\sum_{j=1}^\infty
\left(
\frac
{
\capinRn(K_{j})
}
{
\capinRn(B_{j})
}
\right)^\metricexp
\nonumber
\\
&
&
+
\sum_{j=1}^\infty
r_j
\left(
\sum_{l=1}^j
\frac{\capinRn (K_l)}{r_l^{n-2}}
\right)^\metricexp
\nonumber
\\
&
=
&
I+II.
\end{eqnarray}
Thus to prove 
(\ref{Sigma1})
we need to bound
$I$
and 
$II$.

{\bf 7.}
Utilising
(\ref{sertointegr})
we deduce at once that
\begin{equation*}
I
\lesssim
\int_0^1
\left(
\frac
{
\capinRn(K\cap B(0,r))
}
{
\capinRn(B(0,r))
}
\right)
^
{\metricexp}
\,
\frac{dr}{r}.
\end{equation*}
To estimate
$II$
we
define the function
$\Phi\colon (0,1)\to {\bf R}^1$
by writing
\begin{equation*}
\Phi(r)
=
\capinRn
(K\cap B(0,r)),
\quad
0<r<1.
\end{equation*}
First assume that
$n\geq 4$
and hence
\begin{equation*}
\frac{2}{n-2}\leq 1.
\end{equation*}
In this case 
by the simple change
of the summation order
we discover that
\begin{eqnarray*}
II
&
\leq
&
\sum_{j=1}^\infty
r_j
\,
\sum_{l=1}^j
\left(
\frac{\capinRn (K_l)}{r_l^{n-2}}
\right)
^\metricexp
\\
&
\lesssim
&
\sum_{l=1}^\infty
\left(
\frac{\capinRn (K_l)}{r_l^{n-2}}
\right)
^\metricexp
r_l
\\
&
\lesssim
&
\int_0^1
\left(
\frac
{
\Phi(r)
}
{
r^{n-2}
}
\right)
^\metricexp
\, 
dr
.
\end{eqnarray*}
Assume next that
$n=3$,
and hence
\begin{equation*}
\frac{2}{n-2}=2
\end{equation*}
Then Hardy's inequality
implies
that
\begin{eqnarray*}
II
&
\lesssim
&
\int_0^1
\left(
\int_t^1
\frac{\Phi(r)}{r}
\,
\frac{dr}{r}
\right)^2
\,
dt
\\
&
\lesssim
&
\int_0^1
\left(
\frac{\Phi(r)}{r}
\right)^2
\, dr.
\end{eqnarray*}
Thus for any
$n\geq 3$
we have
\begin{equation*}
II
\lesssim
\int_0^1
\left(
\frac
{
\Phi(r)
}
{
r^{(n-2)/2}
}
\right)^\metricexp
\,
\frac{dr}{r}
.
\end{equation*}
Returning to
(\ref{Sigma1proof2})
and recalling
(\ref{capball})
we derive
\begin{equation*}
\int_{B\setminus \Sigma}
u(x)^{2/(n-2)}
\,
\frac{1}{|x|^{n-1}}
\, 
dx
\lesssim
\int_0^1
\left(
\frac
{
\capinRn(K\cap B(0,r))
}
{
\capinRn (B(0,r))
}
\right)^{2/(n-2)}
\,
\frac{dr}{r}.
\end{equation*}
Employing
(\ref{thinKJ})
we establish
(\ref{Sigma1}).
This completes the proof
of the implication
$(i)\Rightarrow (ii)$
in Theorem~\ref{maintheorem}.

%
%

%

\begin{thebibliography}{}

%
%



\bibitem{Adams_Hedberg_book}
D. R. Adams, L. I. Hedberg,
{\it Function spaces and potential theory},
Springer-Verlag, Berlin, Heidelberg, 1996.


\bibitem{Au1}
T. Aubin,
Equations differentielles non 
lineaires et probleme de Yamabe concernant 
la courbure scalaire. J. Math. Pures Appl. 55, 269--296 (1976)





\bibitem{Aubin_book}
T. Aubin,
\textit{Some nonlinear problems in Riemannian geometry},
Springer-Verlag, Berlin, 1998.




\bibitem{Aviles_CPDE}
P. Aviles,  A study of the singularities of solutions of a class of 
nonlinear elliptic partial differential equations.  
Comm. Partial Differential Equations  7  (1982), no. 6, 609--643.



\bibitem{Aviles_McOwen_Duke_88}
P. Aviles, R. McOwen, 
Complete conformal metrics with negative 
scalar curvature in compact Riemannian manifolds.  Duke Math. J.  56  (1988),  no. 2, 395--398. 



\bibitem{Aviles_McOwen_JDG_85}
P. Aviles, R. McOwen, 
Conformal deformations of complete manifolds with 
negative curvature.  J. Differential Geom.  21  
(1985),  no. 2, 269--281.


\bibitem{Aviles_McOwen_JDG_88}
P. Aviles, R. McOwen, 
Conformal deformation to constant 
negative scalar curvature on noncompact Riemannian 
manifolds.  J. Differential Geom.  27  (1988),  
no. 2, 225--239.





\bibitem{Delanoe_Contemp_Math_90}
Delano\"e, Philippe 
Generalized stereographic projections with prescribed scalar curvature.  
Geometry and nonlinear partial differential equations (Fayetteville, AR, 1990),  
17--25, Contemp. Math., 127, Amer. Math. Soc., Providence, RI, 1992. 


\bibitem{DLeG}
J.-S. Dhersin,  J.-F. Le Gall, Wiener's test for
super-Brownian motion and the Brownian snake, 
Probab. Theory Related Fields {\bf 108} (1997),
103--129.

\bibitem{Dynkin_book}
Dynkin, E. B. 
Diffusions, superdiffusions and partial differential equations. 
American Mathematical Society Colloquium Publications, 50. 
American Mathematical Society, Providence, RI, 2002.



\bibitem{Dynkin_process_manifolds}
Dynkin, E. B. Kuznetsov, E.
Extinction of superdiffusions and semilinear partial differential equations. 
J. Funct. Anal. 162 (1999), no. 2, 346--378.

\bibitem{Escobar_92_1}
J. F. Escobar,
The Yamabe problem on manifolds with boundary.  J. Differential Geom.  35  (1992),  no. 1, 21--84.



\bibitem{Escobar_92_2}
J. F. Escobar,
Conformal deformation of a Riemannian metric to a scalar flat metric with constant mean curvature on the boundary.  Ann. of Math. (2)  136  (1992),  no. 1, 1--50.

\bibitem{Escobar_03}
J. F. Escobar,
Uniqueness and non-uniqueness of metrics with prescribed scalar and mean curvature on compact manifolds with boundary.  J. Funct. Anal.  202  (2003),  no. 2, 424--442.


\bibitem{Finn_94}
D. L. Finn,
Positive solutions of $\Delta\sb g u=u\sp q+Su$ singular at 
submanifolds with boundary.  Indiana Univ. Math. J.  43  (1994),  no. 4, 1359--1397.


\bibitem{Finn_99}
D. L. Finn,
On the negative case of the singular 
Yamabe problem.  J. Geom. Anal.  9  
(1999),  no. 1, 73--92.


\bibitem{Finn_00}
D. L. Finn,
Behavior of positive solutions to 
$\Delta\sb gu=u\sp q+Su$
with prescribed singularities, 
Indiana Univ. Math. J. \textbf{49} (2000),  177--219.


\bibitem{Finn_McOwen_IUMJ_93}
 Finn, David L.; 
 McOwen, Robert C. Singularities and asymptotics for the equation $\Delta\sb gu-u\sp q=Su$.  
 Indiana Univ. Math. J.  42  (1993),  no. 4, 1487--1523.


\bibitem{Gamara_Yacoub}
 N. Gamara, 
 R. Yacoub, 
 CR Yamabe conjecture---the conformally flat case.  
 Pacific J. Math.  201  (2001),  no. 1, 121--175. 



\bibitem{Grigor'yan}
A. Grigor'yan,  
Analytic and geometric background of recurrence and 
non-explosion of the Brownian motion on Riemannian manifolds.  
Bull. Amer. Math. Soc. (N.S.)  36  (1999),  no. 2, 135--249.





\bibitem{Jerison_Lee},
D. Jerison, J. M.  Lee,  
The Yamabe problem on CR manifolds.  
J. Differential Geom.  25  (1987),  no. 2, 167--197.


\bibitem{Kato_Nayatani} 
S. Kato, S. Nayatani,  
Complete conformal metrics with prescribed scalar 
curvature on subdomains of a compact manifold.  Nagoya Math. J.  132  (1993), 155--173.


\bibitem{Kell}
Keller, J. B., 
On solutions of $\Delta u = f(u)$,  
Comm. Pure Appl. Math. 10 (1957), 503--510.



\bibitem{KMPS}
N. Korevaar, R. Mazzeo, 
F. Pacard, R. Schoen,
Refined asymptotics for constant scalar 
curvature metrics with isolated singularities.  
Invent. Math.  135  (1999),  no. 2, 233--272.
 


\bibitem{Koskela} 
P. Koskela, Old and new on the quasihyperbolic metric.  
Quasiconformal mappings and analysis (Ann Arbor, MI, 1995),  
205--219, Springer, New York, 1998.



\bibitem{Lab2}
D. A. Labutin,
Wiener regularity for large solutions of 
nonlinear equations.  Ark. Mat.  41  (2003),  
no. 2, 307--339.


\bibitem{Labflat}
D. A. Labutin,
Polar sets for scalar flat
complete metrics,
in preparation.


\bibitem{Lee_Parker}
J. M. Lee,
T. H. Parker,
The Yamabe problem,
Bull. Amer. Math. Soc. (N.S.) \textbf{17} (1987),  37--91.


\bibitem{LeGall_book}
J.-F. LeGall,
{\it Spatial branching processes, random snakes and
partial differential equations}, Lectures in 
Mathematics ETH Z\"urich. Birkh\"auser Verlag, Basel,
1999. 


\bibitem{LeGcongress}
 Le Gall, Jean-Franois Branching processes, random 
trees and superprocesses. Proceedings of the International 
Congress of Mathematicians, Vol. III (Berlin, 1998).  Doc. Math.  1998,  Extra Vol. III, 279--289



\bibitem{Loewner_Nirenberg}
C. Loewner, L. Nirenberg,
Partial differential equations invariant
under conformal or projective transformations, 
Contributions to analysis (a collection of papers
dedicated to Lipman Bers),  
Academic Press, New York, 1974,
245--272.
  
 \bibitem{Ma_McOwen} 
 Ma, Xiaoyun; McOwen, Robert C. 
 Complete conformal metrics with zero scalar curvature.  
 Proc. Amer. Math. Soc.  115  (1992),  no. 1, 69--77.


\bibitem{Marcus_Veron_Poincare_journal}
Marcus, M. and V\'eron, L., Uniqueness and asymptotic behavior 
of solutions with boundary blow-up for a class of nonlinear 
elliptic equations,  Ann. Inst. H. Poincar\'e Anal. Non Lin\'eaire 14 (1997), 237--274.


\bibitem{Maz'ya_book}
V. G. Maz'ja
{\it Sobolev spaces},
Springer-Verlag, Berlin-New York, 1985.



\bibitem{Mazzeo_IUMJ_91}
R. Mazzeo,  
Regularity for the singular 
Yamabe problem.  Indiana Univ. Math. J.  40  (1991),  
no. 4, 1277--1299.


\bibitem{Mazzeo_Pacard_Duke}
R. Mazzeo, F. Pacard,
Constant scalar curvature metrics with
isolated singularities, 
Duke Math. J. 
{\bf 99} (1999), 353--418.


\bibitem{Mazzeo_Pacard_JDG}
R. Mazzeo, F. Pacard,
A construction of singular solutions 
for a semilinear elliptic equation using 
asymptotic analysis.  J. Differential Geom.  44  
(1996),  no. 2, 331--370.


\bibitem{Mazzeo_Pollack_Uhlenbeck_JAMS}
R. Mazzeo, D. Pollack, K. Uhlenbeck,
Moduli spaces of
singular Yamabe metrics, 
J. Amer. Math. Soc. {\bf 9} (1996), no. 2, 
303--344.

\bibitem{Mazzeo_Pollack_Uhlenbeck_TPMNA}
R. Mazzeo, D. Pollack, K. Uhlenbeck,
Connected sum constructions for constant scalar curvature metrics.  
Topol. Methods Nonlinear Anal.  6  (1995),  no. 2, 207--233. 




\bibitem{Mazzeo_Smale_JDG}
Mazzeo, Rafe; Smale, Nathan 
Conformally flat metrics of constant 
positive scalar curvature on subdomains of the sphere.  
J. Differential Geom.  34  (1991),  no. 3, 581--621.



\bibitem{Mazzeo_Taylor_uniformisation}
Mazzeo, Rafe; Taylor, Michael 
Curvature and uniformization.  Israel J. Math.  130  (2002), 323--346.




\bibitem{McOwen_survey}
R. C. McOwen,
Results and open questions on the singular
Yamabe problem,
Dynamical systems and differential equations, Vol. II 
(Springfield, MO, 1996). 
Discrete Contin. Dynam. Systems 1998, 
Added Volume II, 123--132.

\bibitem{Moser_1}
J. Moser,  
On Harnack's theorem for elliptic differential equations.  Comm. Pure Appl. Math.  14  1961 577--591.


\bibitem{Moser_2}
J. Moser
A new proof of De Giorgi's theorem concerning the regularity 
problem for elliptic differential equations.  Comm. Pure Appl. Math.  13  1960 457--468.
 
 \bibitem{Osserm}
Osserman, R., On the inequality $\Delta u \geq f(u)$,  
Pacific J. Math. 7 (1957), 1641--1647.



\bibitem{Pacard_TMNA}
Pacard, Frank 
The Yamabe problem on subdomains of even-dimensional spheres.  
Topol. Methods Nonlinear Anal.  6  (1995),  no. 1, 137--150.

\bibitem{Schoen_Yamabe}
R. Schoen,
Conformal deformation of a Riemannian metric to
constant scalar curvature, 
J. Differential Geom. 
\textbf{20} (1984), 479--495.


\bibitem{Schoen_ICM}
R. Schoen,
Recent progress in geometric partial differential equations.  
Proceedings of the International Congress of Mathematicians, Vol. 1, 2 (Berkeley, Calif., 1986),  121--130, Amer. Math. Soc., Providence, RI, 1987.


\bibitem{Schoen_CPAM}
R. Schoen,
The existence of weak solutions with prescribed
singular behavior for a conformally invariant 
scalar equation, Comm. Pure Appl. Math. 
\textbf{41}
(1988),  317--392. 



\bibitem{Schoen_Yau_Invent_88}
R. Schoen, S. T. Yau,
Conformally flat manifolds, Kleinian groups and
scalar curvature, 
Invent. Math.  
\textbf{92} (1988), 47--71. 




\bibitem{Schoen_Yau_book}
R. Schoen, S. T. Yau,
{\it Lectures on differential geometry},
Conference Proceedings and Lecture Notes in Geometry
and Topology, I. 
International Press, Cambridge, MA, 1994. 

\bibitem{Trudinger_weak_Harnack}
 Trudinger, Neil S. On Harnack type inequalities 
and their application to quasilinear elliptic equations.  
Comm. Pure Appl. Math.  20  1967 721--747.


\bibitem{Tru}
N. S. Trudinger,
Remarks concerning the conformal deformation 
of Riemannian structures on compact
manifolds, 
Ann. Scuola Norm. Sup. Pisa (3)
\textbf{22} (1968), 265--274. 


\bibitem{Veron_JDE}
 L. V\'eron,  
Singularites eliminables d'equations elliptiques non 
lineaires. 
(French)  
J. Differential Equations  41  (1981), no. 1, 87--95.



\bibitem{Wiener}
N. Wiener,
Collected works. Vol. I.
Mathematical philosophy and foundations; potential theory; Brownian movement, Wiener integrals, ergodic and chaos theories, turbulence and statistical mechanics. 
MIT Press, Cambridge, Mass.-London, 1976.



\bibitem{Yam}
H. Yamabe,
On a deformation of Riemannian
structures on compact manifolds,
Osaka  Math. J. \textbf{12} (1960), 
21--37.


\bibitem{Yau_problems}
S. T. Yau,
Open problems in geometry.
Differential geometry: partial differential equations on manifolds
Proceedings of Symp. Pure Math.,
{\bf 54} Part 1, 1--28,
Amer. Math. Soc., Providence, RI, 1993.


\bibitem{Ziemer_book}
W. P. Ziemer, 
{\it
Weakly differentiable functions. Sobolev spaces 
and functions of bounded variation.} 
Graduate Texts in Mathematics, 120. Springer-Verlag, New York, 1989.



\end{thebibliography}
%

\end{document}